\documentclass[12pt, reqno]{amsart}

\usepackage[margin=1in,showframe=false]{geometry}

\usepackage{amsmath, amsthm, amsfonts, amssymb, graphicx, color, bm}
\usepackage{mathtools}
\usepackage{microtype}
\usepackage{float}
\usepackage[table]{xcolor}
\usepackage{enumitem}
\usepackage{algorithm}
\usepackage{algorithmic}
\usepackage[bookmarksnumbered, colorlinks, plainpages]{hyperref}
\usepackage{tikz}
\usetikzlibrary{calc,angles,quotes}

\newtheorem{theorem}{Theorem}[section]

\theoremstyle{definition}

\theoremstyle{remark}

\numberwithin{equation}{section}
\numberwithin{figure}{section}

\begin{document}

\title[Reconstruction by diffracted waves]{Numerical reconstruction of polygonal obstacles\\ by diffracted waves}

\author[X. Chen]{Xi Chen}
\address{Shanghai Center for Mathematical Sciences, Fudan University, Shanghai 200438, China; %
    School of Mathematical Sciences, Fudan University, Shanghai 200433, China; %
	Center for Applied Mathematics, Fudan University, Shanghai 200433, China.}
\email{xi\_chen@fudan.edu.cn}

\author[K. Xiong]{Kuangmiao Xiong}
\address{Shanghai Center for Mathematical Sciences, Fudan University, Shanghai 200438, China.}
\email{24110840017@m.fudan.edu.cn}

\author[B. Xu]{Boxi Xu}
\address{School of Mathematics, Shanghai University of Finance and Economics, Shanghai 200433, China.}
\email{xu.boxi@mail.sufe.edu.cn}

\begin{abstract}
We present a numerical method for recovering rigid polygonal obstacles in two dimensions, based on the diffraction-based inverse scattering theory of Bao, Chen, Lu, and Xiong. %
A single pulse is emitted from a point source, and the resulting wave field is measured on a receiver curve. %
After isolating the diffracted signals, we recover the visible vertices from arrival times and the corner orientations and opening angles from amplitudes. %
The use of amplitude information permits local reconstruction in the reflection shadow, where reflected waves do not reach the receiver. %
The amplitudes also determine the direction of a shadowed edge adjacent to a visible vertex, enabling local non-line-of-sight imaging. %
Numerical experiments for a convex corner, a concave corner, a triangle, a pentagon, and a star-shaped obstacle confirm that the visible portion of the boundary can be recovered stably from diffracted signals. %
A further convex-corner test with additive noise at level $0.2\%$ on the single arc $L_{1}$ confirms that the tip and opening remain recoverable. %
\end{abstract}



\maketitle


\section{Introduction}
\label{sec:introduction}

\subsection{Overview}

Inverse obstacle scattering aims to recover the location and shape of unknown obstacles from wave measurements. %
Classical inversion approaches \cite{CK-book, KKL-book} usually rely on far-field patterns or on boundary measurements of full wave fields to ensure uniqueness and stability. %
Bao, Chen, Lu, and Xiong \cite{BCLX1} developed an alternative framework that exploits the diffracted waves by vertices: the visible portion of a polygonal obstacle in $\mathbb{R}^{2}$ is uniquely determined by the vertex-diffracted waves from a single pulse recorded on a receiver curve of essentially arbitrary size and placement. %
Two features distinguish this diffraction imaging from reflection tomography. %
First, tip-diffracted waves propagate omnidirectionally from a vertex, so they can be captured by an arbitrarily small receiver whenever that vertex is visible. %
Consequently, the required data is minimal compared with conventional wave measurements. %
Second, the uniqueness proof is constructive: the arrival times of diffracted signals determine the vertex locations, and their amplitudes determine the orientations and opening angles of the corners, with both quantities obtainable via Cheeger–Taylor functional calculus on cones. %
The present paper implements this constructive procedure numerically. %
We generate synthetic data by a staggered-grid finite-difference scheme with PML and Neumann boundary treatment. %
The inversion proceeds by extracting the diffracted wave data from the full field and then solving two successive nonlinear least-squares problems. %


\subsection{Inverse obstacle scattering by diffracted waves}

Consider the following initial-boundary value problem for the wave equation. %
Let the open set $K \subset \mathbb{R}^{2}$ be a fixed rigid polygonal obstacle and denote by $K^{\mathrm{ext}} := \mathbb{R}^{2} \setminus \overline{K}$ its exterior domain. %
The scattering of waves from a point source $\delta_{O}$ by the obstacle $K$ is described by %
\begin{equation}
\label{eqn:model}
\left\{~
\begin{aligned}
  ( \partial_{tt} - \Delta ) u(t,x;O) &= 0,
    & ~ & (t,x) \in \mathbb{R} \times K^{\mathrm{ext}}, \\
  \partial_{\nu} u(t,x;O) &= 0,
    & ~ & (t,x) \in \mathbb{R} \times \partial K, \\
  u(0,x;O) &= 0,
    & ~ & x \in K^{\mathrm{ext}}, \\
  \partial_{t} u(0,x;O) &= \delta_{O},
    & ~ & x \in K^{\mathrm{ext}}, \\
\end{aligned}
\right.
\end{equation}
where $u(t,x;O)$ denotes the exterior wave field, %
$O \in K^{\mathrm{ext}}$ is the source point, %
and $\nu$ is the unit outward normal on $\partial K$. %

The inverse problem is to recover the visible portion of $K$ from the wave field $u(t,x;O)$ recorded on a receiver curve $L \subset K^{\mathrm{ext}}$. %
A vertex $P_{i}$ of $K$ is called \textit{visible} with respect to $O$ and $L$ if the line segments joining $P_{i}$ to $O$ and to any segment of $L$ lie entirely in $K^{\mathrm{ext}}$.

\begin{figure}[htb]
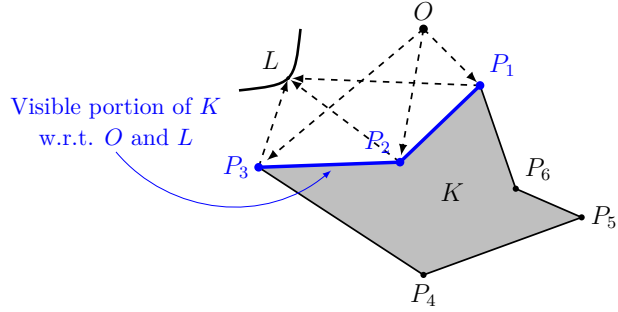

\centering
\resizebox{0.5\textwidth}{!}{
\tikz[>=latex]{
\coordinate (O)  at (0,0);
\coordinate (P1) at (-45:1.3);
\coordinate (P2) at (-100:2.2);
\coordinate (P3) at (-140:3.5);
\coordinate (P4) at (-90:4);
\coordinate (P5) at (-50:4);
\coordinate (P6) at (-60:3);
\coordinate (K)  at (-80:2.7);
\coordinate (L)  at (-2.2,-0.8);
\coordinate (F)  at (-5,-1.5);
%
\draw[thick,fill=lightgray] (P1) -- (P2) -- (P3) -- (P4) -- (P5) -- (P6) -- cycle;
\draw[very thick] (-2,0) .. controls (-2.1,-0.9) .. (-3,-1);
\draw[ultra thick,blue] (P1) -- (P2) -- (P3);
\draw[thick,  red!0!black,dashed,->] (O)  -- ($ (P1) !0.05! (O) $);
\draw[thick,  red!0!black,dashed,->] (P1) -- ($ (L) !0.02! (P1) $);
\draw[thick, blue!0!black,dashed,->] (O)  -- ($ (P2) !0.05! (O) $);
\draw[thick, blue!0!black,dashed,->] (P2) -- ($ (L) !0.02! (P2) $);
\draw[thick,green!0!black,dashed,->] (O)  -- ($ (P3) !0.05! (O) $);
\draw[thick,green!0!black,dashed,->] (P3) -- ($ (L) !0.02! (P3) $);
\fill[black] (O)  circle (0.07) node[above] {$O$};
\fill[ blue] (P1) circle (0.07) node[above right] {$P_{1}$};
\fill[ blue] (P2) circle (0.07) node[above left] {$P_{2}$};
\fill[ blue] (P3) circle (0.07) node[left] {$P_{3}$};
\fill[black] (P4) circle (0.05) node[below] {$P_{4}$};
\fill[black] (P5) circle (0.05) node[right] {$P_{5}$};
\fill[black] (P6) circle (0.05) node[above right] {$P_{6}$};
\node at (K) {$K$};
\fill[black] (L)  circle (0.03) node[above left] {$L$};
\node[align=center,blue] at (F) {\small Visible portion of $K$ \\ \small w.r.t. $O$ and $L$};
\draw[blue,->] (F)++(0,-0.5) to[bend right=45] (-1.5,-2.3);
}}
\caption{The vertices $P_{1}$, $P_{2}$ and $P_{3}$ are visible with respect to $O$ and $L$.}
\label{fig:visible}
\end{figure}


\subsubsection{Corner identification}

\begin{figure}[htb]
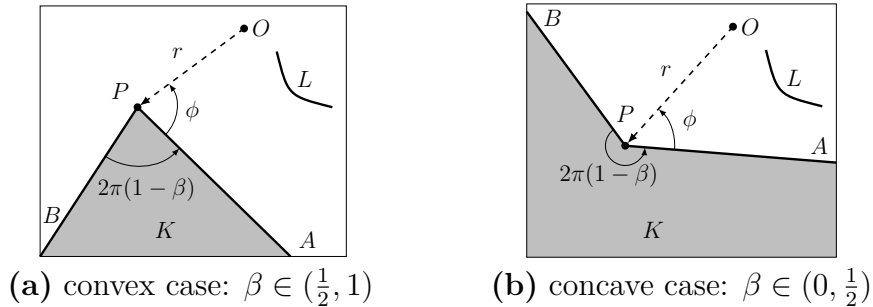

\centering
\begin{tabular}{ccc}
\resizebox{0.25\textwidth}{!}{
\tikz[>=latex]{
\coordinate (O)  at (45:2);
\coordinate (P)  at (-0.5,0);
\coordinate (P1) at (-50:3.5);
\coordinate (P2) at (-130:3.5);
\coordinate (K)  at (-90:2.2);
\coordinate (L)  at (2.5,0.5);
\fill[lightgray] (P1) -- (P) -- (P2);
\draw[very thick] (P1) -- (P) -- (P2);
\draw[very thick] (2,1) .. controls (2.2,0.2) .. (3,0);
\draw[thick,dashed,->] (O) -- ($ (P) !0.02! (O) $) node[pos=0.5,above left] {$r$};
\draw pic["$\phi$",draw,->,angle radius=20,angle eccentricity=1.4] {angle=P1--P--O};
\draw pic["$2\pi(1-\beta)$",draw,->,angle radius=30,angle eccentricity=1.4] {angle=P2--P--P1};
\fill[black] (O) circle (0.07) node[right] {$O$};
\fill[black] (P) circle (0.07) node[above left] {$P$};
\node at (L) {$L$};
\node at (K) {$K$};
\node[above right] at (P1) {$A$};
\node[xshift=6,yshift=21] at (P2) {$B$};
\draw (P2) -- ($ (P2) + (5.5,0) $) -- ++(0,4.5) -- ($ (P2) + (0,4.5) $) -- cycle;
}}
& &
\resizebox{0.25\textwidth}{!}{\tikz[>=latex]{
\coordinate (O)  at (45:2);
\coordinate (P)  at (-0.5,-0.7);
\coordinate (P1) at ($ (0,-1) + (0:3.25) $);
\coordinate (P2) at ($ (0,-1) + (130:3.5) $);
\coordinate (P3) at (-130:3.5);
\coordinate (K)  at (-90:2.2);
\coordinate (L)  at (2.5,0.5);
\fill[lightgray] (P1) -- (P) -- (P2) -- (P3) -- ($ (P3) + (5.5,0) $) -- cycle;
\draw[very thick] (P1) -- (P) -- (P2);
\draw[very thick] (2,1) .. controls (2.2,0.2) .. (3,0);
\draw[thick,dashed,->] (O) -- ($ (P) !0.02! (O) $) node[pos=0.5,above left] {$r$};
\draw pic["$\phi$",draw,->,angle radius=25,angle eccentricity=1.4] {angle=P1--P--O};
\draw pic["$2\pi(1-\beta)$",draw,->,angle radius=10,angle eccentricity=1.7] {angle=P2--P--P1};
\fill[black] (O) circle (0.07) node[right] {$O$};
\fill[black] (P) circle (0.07) node[above=8] {$P$};
\node at (L) {$L$};
\node at (K) {$K$};
\node[above left] at (P1) {$A$};
\node[xshift=13,yshift=-3] at (P2) {$B$};
\draw (P3) -- ($ (P3) + (5.5,0) $) -- ++(0,4.5) -- ($ (P3) + (0,4.5) $) -- cycle;
}}
\\
\textbf{(a)} convex case: $\beta \in (\frac{1}{2},1)$
& \qquad\qquad &
\textbf{(b)} concave case: $\beta \in (0,\frac{1}{2})$ \\
\end{tabular}
\caption{The obstacle $K$ with a single corner.}
\label{fig:corner}
\end{figure}

We first consider a single corner with tip $P$ (Figure \ref{fig:corner}). %
In polar coordinates centred at $P$, with the polar axis along one edge of the corner, we write $O := (r,\phi)$ with $r := |O-P|$, where $\phi$ is the angle between that edge and the ray $PO$, measured through the exterior sector. %
Let $\beta \in (0,1]$ be the opening parameter such that the exterior angle at $P$ equals $2\pi \beta$ and the interior angle equals $2\pi (1-\beta)$. %
Thus, in Figure \ref{fig:corner}, $K$ is called \textit{convex} if $\beta \in (\tfrac{1}{2},1)$ (interior angle in $(0,\pi)$) and \textit{concave} if $\beta \in (0,\tfrac{1}{2})$ (interior angle in $(\pi,2\pi)$). %
The corner is called \textit{non-diffractive} if $\beta \in \{ N^{-1} : N \in \mathbb{Z}^{+} \}$ and \textit{diffractive} otherwise.%
\footnote{The diffracted waves vanish when the corner is non-diffractive; see \cite[Remarks 1 and 7]{BCLX1}.} %

Measurements of the principal part of the diffracted wave on a receiver curve $L$ determine a diffractive corner uniquely, as stated next.

\begin{theorem}[\cite{BCLX1}]
\label{thm:uniqueness}
Let $K \subset \mathbb{R}^{2}$ be a diffractive corner with tip $P$, let $O \in K^{\mathrm{ext}}$, and let $L \subset K^{\mathrm{ext}}$ be a receiver curve such that $O$ and $L$ do not lie on the two branches of a hyperbola.%
\footnote{If $O$ and $L$ lie on the two branches of a hyperbola, two distinct tips may give the same diffracted arrival time on $L$; see \cite[Remark 2]{BCLX1}.} %
Write $r = |O-P|$ and $r_{1}(x) := |x-P|$. %
Then the mapping %
\begin{equation}
\label{eqn:uniqueness}
  (P,\phi,\beta) ~\longmapsto~ \bigl\{ u_{\mathcal{D}}(t,x;O) \,:\, t = r + r_{1}(x), ~ x \in L \bigr\}
\end{equation}
is injective, where $u_{\mathcal{D}}$ is the principal part of the diffracted wave, given in \eqref{eqn:diffract}. %
\end{theorem}


\subsubsection{Polygon identification}

For a polygonal obstacle, these diffraction data determine each visible corner under suitable geometric assumptions. %

\begin{theorem}[\cite{BCLX1}]
\label{thm:polygon}
Let $K \subset \mathbb{R}^{2}$ be a compact convex polygon, let $O \in K^{\mathrm{ext}}$ emit no grazing ray,%
\footnote{A ray issuing from $O$ is called a \textit{grazing ray} if the line through $O$ is tangent to $\partial K$. %
Along a grazing ray, the wave field has the same Sobolev regularity as the incident wave, so primary and multiple diffracted waves cannot be distinguished by their regularity. %
We therefore exclude grazing rays.} %
and let $L \subset K^{\mathrm{ext}}$ be a receiver curve that is not assembled solely from straight segments and hyperbolic arcs. %
Denote by $\{ P_{j} \}_{j=1}^{N}$ the visible vertices of $K$ with respect to $O$ and $L$, and write $(P_{j},\phi_{j},\beta_{j})$ for the corresponding polar parameters as in Theorem \ref{thm:uniqueness}. %
Then the mapping sending $\{ (P_{j},\phi_{j},\beta_{j}) \}_{j=1}^{N}$ to the localised wave field on $L$, observed at the successive diffracted arrival times, is injective. %
\end{theorem}

The proof of Theorem \ref{thm:polygon} is constructive. %
On a suitable open sub-arc of $L$ the primary diffracted waves from distinct vertices are pairwise separated in time and arrive after geometric waves. %
These primary diffracted waves can then be distinguished from multiple diffracted by their Sobolev regularities. %
For each separated primary diffracted wave, an analysis based on finite propagation speed yields the same principal parts as in the single-corner case, allowing reconstruction by Theorem \ref{thm:uniqueness}. %
Assembling the recovered wedges reconstructs precisely the visible portion of $\partial K$. %

\begin{figure}[htb]
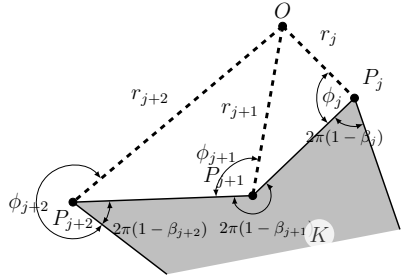

\centering
\resizebox{0.32\textwidth}{!}{
\tikz[>=latex,scale=1.5]{
\coordinate (O)  at (0,0);
\coordinate (P1) at (-45:1.3);
\coordinate (P2) at (-100:2.2);
\coordinate (P3) at (-140:3.5);
\coordinate (P4) at (-115:3.5);
\coordinate (P6) at (-60:3);
\coordinate (K)  at (-70:2.7);
%
\draw[thick,fill=lightgray] (P6) -- (P1) -- (P2) -- (P3) -- (P4);
\draw[ultra thick,black,dashed] (O) -- (P1) node[pos=0.4,above right] {$r_{j}$};
\draw pic["$\phi_{j}$",draw,<->,angle radius=20,angle eccentricity=0.6] {angle=O--P1--P2};
\draw pic["{\scriptsize$2\pi(1-\beta_{j})$}",draw,<->,angle radius=15,angle eccentricity=1.5] {angle=P2--P1--P6};
\draw[ultra thick,black,dashed] (O) -- (P2) node[pos=0.5,left] {$r_{j+1}$};
\draw pic["$\phi_{j+1}$",draw,<->,angle radius=20,angle eccentricity=1.5] {angle=O--P2--P3};
\draw pic["{\scriptsize$2\pi(1-\beta_{j+1})$}",draw,<->,angle radius=10,angle eccentricity=2] {angle=P3--P2--P1};
\draw[ultra thick,black,dashed] (O) -- (P3) node[pos=0.5,above left] {$r_{j+2}$};
\draw pic["$\phi_{j+2}$",draw,<->,angle radius=20,angle eccentricity=1.2] {angle=O--P3--P4};
\draw pic["{\scriptsize$2\pi(1-\beta_{j+2})$}",draw,<->,angle radius=20,angle eccentricity=2.5] {angle=P4--P3--P2};
%
\draw[fill] (O)  circle (0.05) node[above] {$O$};
\draw[fill] (P1) circle (0.05) node[above right] {$P_{j}$};
\draw[fill] (P2) circle (0.05) node[above left] {$P_{j+1}$};
\draw[fill] (P3) circle (0.05) node[below] {$P_{j+2}$};
\node at (K) {$K$};
}}
\caption{Determining the location and shape of the visible portion of $K$ in $\mathbb{R}^2$}
\label{fig:polar}
\end{figure}

We remark that convexity in Theorem \ref{thm:polygon} is assumed for brevity of the argument. %
The local single-corner analysis remains valid at every diffractive visible vertex of a non-convex polygon, provided $\beta_{j} \in (0,1] \setminus \{ N^{-1} : N \in \mathbb{Z}^{+} \}$ and grazing rays are avoided. %
This is confirmed numerically for a star-shaped obstacle in Section \ref{sec:numerical}. %

The constructive character of Theorems \ref{thm:uniqueness} and \ref{thm:polygon} is essential for computation: the wave field on $L$ consists of incident, reflected and diffracted waves, whereas the recovery uses only the primary diffracted waves. %
Incident fields carry no information about $K$, and reflected fields may be missed by a small receiver. %
Finite propagation speed separates the diffracted and reflected waves from distinct vertices and edges, and their Sobolev regularities distinguish geometric waves, primary diffracted waves and multiple diffracted waves. %
These principles permit the use of receiver curves of essentially arbitrary size and placement and form the basis of the reconstruction scheme in Section \ref{sec:reconstruction-scheme}. %


\subsection{Related work}

Inverse problems for polygonal obstacles have been extensively studied within classical inverse scattering theory. %
With the far-field pattern of a single incident wave, \cite{CY-IP-2003, EY-IP-2006, HL-IP-2014} established unique reconstruction of impenetrable polygonal and polyhedral obstacles. %
These uniqueness results for a single incident wave were later extended to penetrable scatterers \cite{HSV-SIMA-2016, EH-IP-2019, HLZ-IP-2019, HV-IP-2020, CX-CPDE-2021, DFLY-IPI-2024}. %
These approaches reconstruct geometry primarily from reflected or transmitted information and typically require far-field or boundary data of substantial aperture. %

The two-dimensional framework in \cite{BCLX1} recovers the visible corners from diffracted waves on an exterior receiver curve. %
In $\mathbb{R}^{3}$, the approach in \cite{BCLX2} reconstructs the visible part of strictly convex smooth cones using diffraction coefficients together with arrival times or directions measured on a small receiver surface patch in the reflection shadow. %

In non-destructive testing, the method of time-of-flight diffraction (TOFD) uses diffracted arrival times to locate and size cracks \cite{SilkLidington1975}. %
Further developments in TOFD include Sparse-SAFT for improved resolution \cite{SunEtAl2022}, and a forward surrogate model for signals from rough defects \cite{PaialungaEtAl2024}. %

Regarding the forward theory, the mathematical analysis of diffraction originates with Sommerfeld \cite{So96, So54} and Friedlander \cite{Fr58}. %
Keller's geometrical theory of diffraction \cite{Ke62} provides the classical asymptotic picture. %
Cheeger--Taylor \cite{CT82a, CT82b} computed the wave kernel on product-type cones by functional calculus, from which the principal part $u_{\mathcal{D}}$ is obtained. %
On the computational side we rely on staggered-grid leapfrog schemes \cite{Virieux1986}, split-field PML \cite{Berenger1994, Collino2001}, and ghost-point Neumann treatments \cite{Gibou2002, Coco2013, Arias2018}. %
The resulting forward solver supplies synthetic data to a diffraction inversion rather than to a full-waveform inversion. %


\subsection{Contributions}

The contributions of the paper are threefold. %
First, we implement the constructive theory in Theorems \ref{thm:uniqueness} and \ref{thm:polygon} as a two-step least-squares scheme: arrival times locate visible vertices, and amplitudes recover the angular parameters $(\phi,\beta)$ of each corner. %
The inversion uses explicit diffraction formulas and requires no wave-equation solves. %
This permits local non-line-of-sight imaging in the reflection shadow, provided the vertex is visible with respect to the source and receiver curve. %
Second, we verify that the staggered-grid solver reproduces the geometric and diffracted wavefronts predicted by Cheeger--Taylor theory. %
We compare the signals recorded on a receiver curve in Regions A, B, and C, including regions where reflected waves are absent but diffracted waves remain. %
Third, we apply the reconstruction scheme to a convex corner, a concave corner, a triangle, a pentagon and a non-convex star-shaped obstacle, using only diffracted signals recorded on receiver curves. %
For convex corner, quantitative reconstruction errors are reported. %
We also repeat the convex-corner inversion on the single arc after adding noise at level $0.2\%$ to the diffracted records. %


\subsection{Structure of the paper}

Section \ref{sec:reconstruction-scheme} presents the diffraction formula and reconstruction scheme. %
Section \ref{sec:numerical} describes the forward solver, analyses simulated diffracted signals, and reports the numerical reconstruction experiments. %
Section \ref{sec:conclusions} summarises the results and limitations. %


\section{The reconstruction scheme}
\label{sec:reconstruction-scheme}

This section converts the constructive uniqueness theory into a computational pipeline. %
For a single diffractive corner, the principal part $u_{\mathcal{D}}$ of the diffracted wave, given in \eqref{eqn:diffract}, determines the tip $P$ and the angular parameters $(\phi,\beta)$ (Theorem \ref{thm:uniqueness}). %
The principal part has a Heaviside jump at diffracted arrival time $r + r_{1}(x)$, with its amplitude governed by the angular factor $S(\theta_{1}(x);\phi,\beta)$ in \eqref{eqn:S-factor}. %
The angular uniqueness result in \cite[Theorem 4]{BCLX1} uses the coefficient $S_{\phi,\lambda}$ defined in \eqref{eqn:coefficient}, with receiver angles measured relative to the source direction. %
For a convex polygon, each visible vertex generates a primary diffracted wave. %
Once these waves are separated by their diffracted arrival times and Sobolev regularities, each visible corner is reconstructed separately. %
This gives the two successive least-squares problems in Algorithm \ref{alg:two-step}. %

We retain the notation $K$ for the obstacle, $K^{\mathrm{ext}}$ for its exterior, $O \in K^{\mathrm{ext}}$ for the source, and $L \subset K^{\mathrm{ext}}$ for the receiver curve. %
For a corner with tip $P$, we use polar coordinates centred at $P$, with the polar axis along one edge. 
Write $O = (r,\phi)$ and $x = (r_{1},\theta_{1}) \in K^{\mathrm{ext}}$, where 
\begin{equation*}
  r = |O-P|, \qquad r_{1} = r_{1}(x) = |x-P|,
\end{equation*}
and $\theta_{1} = \theta_{1}(x) \in (0,2\pi \beta)$, where the opening parameter is $\beta \in (0,1]$. %
The exterior and interior angles are $2\pi \beta$ and $2\pi (1-\beta)$. %
The corner is non-diffractive when $\beta \in \{ N^{-1} : N \in \mathbb{Z}^{+} \}$, and diffractive otherwise. %


\subsection{Diffraction by a single corner}

For a single corner with tip $P$, the exterior domain is $K^{\mathrm{ext}} = \mathbb{R}^{+} \times (0,2\pi \beta)$, where $\beta \in (\tfrac{1}{2},1)$ for a convex corner and $\beta \in (0,\tfrac{1}{2})$ for a concave corner. %

The total wave field $u$ can be rewritten as %
\begin{equation*}
  u = u^{\mathrm{geo}} + u^{\mathrm{diff}},
\end{equation*}
and the geometric field $u^{\mathrm{geo}} = u^{\mathrm{inc}} + u^{\mathrm{refl}}$ consists of the incident and reflected waves. %

For fixed $t > 0$, write $\mathcal{G}^{O}_{t}$ for the geometric wavefront, $\mathcal{D}^{O}_{t}$ for the diffracted wavefront, and the forward set $\mathcal{F}^{O}_{t} := \overline{\mathcal{G}^{O}_{t}} \cap \mathcal{D}^{O}_{t}$ for their intersection. %
Denote the diffracted arrival time at $x$ by %
\begin{equation}
\label{eqn:arrival_time}
  t_{P}(x) := r + r_{1}(x) = |O-P| + |x-P|.
\end{equation}
Then %
\begin{equation*}
  \mathcal{D}^{O}_{t} := \bigl\{ x \in K^{\mathrm{ext}} : t = r + r_{1}(x) \bigr\},
\end{equation*}
and the exterior is partitioned into %
\begin{equation*}
\begin{aligned}
  \mathrm{Region\ I}: & \quad
    \bigl\{ x \in K^{\mathrm{ext}} : t < |x-O| \bigr\}, \\
  \mathrm{Region\ II}: & \quad
    \bigl\{ x \in K^{\mathrm{ext}} : |x-O| < t < r + r_{1}(x) \bigr\}, \\
  \mathrm{Region\ III}: & \quad
    \bigl\{ x \in K^{\mathrm{ext}} : t > r + r_{1}(x) \bigr\}.
\end{aligned}
\end{equation*}

Away from the forward set $\mathcal{F}^{O}_{t}$, the principal part $u_{\mathcal{D}}$ of $u^{\mathrm{diff}}$ on $\mathcal{D}^{O}_{t}$ is %
\begin{equation}
\label{eqn:diffract}
  u_{\mathcal{D}}(t,r_{1},\theta_{1};r,\phi,\beta)
    = H( t - r - r_{1} ) \, \biggl( - \frac{S(\theta_{1};\phi,\beta)}{16\pi \beta \sqrt{r\,r_{1}}} \biggr),
\end{equation}
where $H$ is the Heaviside function and %
\begin{equation}
\label{eqn:S-factor}
\begin{aligned}
  S(\theta_{1};\phi,\beta)
  &:= \cot \biggl( \frac{\pi+\theta_{1}-\phi}{4\beta} \biggr)
    + \cot \biggl( \frac{\pi-\theta_{1}+\phi}{4\beta} \biggr) \\
  &\quad
    + \cot \biggl( \frac{\pi-\theta_{1}-\phi}{4\beta} \biggr)
    + \cot \biggl( \frac{\pi+\theta_{1}+\phi}{4\beta} \biggr).
\end{aligned}
\end{equation}
Equivalently, the coefficient %
\begin{equation}
\label{eqn:coefficient}
  S_{\phi,\lambda}(\theta_{1}-\phi) = \frac{S(\theta_{1};\phi,\beta)}{2 \beta}
\end{equation}
is used in \cite[Section 3]{BCLX1} with $\lambda := (2\beta)^{-1}$.
For diffractive corners with $\beta \in (0,1] \setminus \{ N^{-1} : N \in \mathbb{Z}^{+} \}$, the coefficient $S_{\phi,\lambda}$ uniquely determines the angular parameters; see \cite[Theorem 4]{BCLX1}. %
Consequently, measurements of the principal part $u_{\mathcal{D}}$ on a receiver curve $L$ uniquely determine $(P,\phi,\beta)$. %


\subsection{Two-step reconstruction for a single corner}
\label{subsec:two-step}

In \eqref{eqn:diffract}, the diffracted arrival time depends on the tip location, whereas the amplitude also depends on the angular parameters. %
We therefore reconstruct these parameters in two steps. %

\medskip\noindent\textbf{Measurement and separation of diffracted waves.}\par\nopagebreak\smallskip
On the receiver curve $L$, one records the total wave field $u(t,x_{i};O)$ at $x_{i} \in L$ ($i = 1,2,\dots,M$). %
A parallel obstacle-free simulation produces the incident field $u^{\mathrm{inc}}(t,x_{i};O)$. %
Subtracting this incident field gives the scattered field %
\begin{equation}
\label{eqn:extract-diff}
  u^{\mathrm{sc}}(t,x_{i};O) := u(t,x_{i};O) - u^{\mathrm{inc}}(t,x_{i};O),
\end{equation}
which contains reflected and diffracted waves ($u^{\mathrm{sc}} = u^{\mathrm{refl}}+u^{\mathrm{diff}}$). %
We then isolate the primary diffracted signals by restricting the measurements to suitable time intervals and, when needed, by a filter based on Sobolev regularity. %
The filter retains the Heaviside singularity $H(t-t_{P}(x))$ of the primary diffracted wave and suppresses the stronger geometric waves. %

The resulting numerical diffracted signal is denoted by $u_{\mathcal{D}}^{\mathrm{num}}$. %
Near each diffracted wavefront, it approximates the principal part $u_{\mathcal{D}}$. %
From each extracted signal $u_{\mathcal{D}}^{\mathrm{num}}$ at $x_{i}$, we estimate the arrival time $T_{i}$ from the onset of the jump. %
For instance, the onset can be detected by thresholding $|\partial_{t} u_{\mathcal{D}}^{\mathrm{num}}|$. %
The signed amplitude $U_{i}$ is then estimated from the jump height across $t = T_{i}$. %
The amplitude approximates the coefficient in \eqref{eqn:diffract}, %
\begin{equation}
\label{eqn:jump}
  U_{i} ~\approx~ -\frac{S(\theta_{1}(x_{i});\phi,\beta)}{16\pi \beta \sqrt{r\,r_{1}(x_{i})}}.
\end{equation}

\medskip\noindent\textbf{Step 1. Reconstruction of the tip location.}\par\nopagebreak\smallskip
With $t_{P}(x) = |O-P| + |x-P|$ as in \eqref{eqn:arrival_time}, solve the nonlinear least-squares problem
\begin{equation}
\label{eqn:step1}
  \min_{P \in \mathbb{R}^{2}} \sum_{i=1}^{M} \bigl| t_{P}(x_{i}) - T_{i} \bigr|^{2}.
\end{equation}
Under the geometric condition in Theorem \ref{thm:uniqueness}, the mapping $P \mapsto t_{P}(x)|_{x \in L}$ is injective. %
Thus \eqref{eqn:step1} recovers the tip uniquely. %

\medskip\noindent\textbf{Step 2. Reconstruction of the angular parameters.}\par\nopagebreak\smallskip
With $P$ recovered, $r = |O-P|$ and $r_{1}(x_{i}) = |x_{i}-P|$ are known. %
In polar coordinates centred at $P$ with the polar axis along the ray $PO$, the relative angles $\alpha_{1}(x_{i}) = \phi - \theta_{1}(x_{i})$ can be computed, although $\phi$ and $\theta_{1}(x_{i})$ remain unknown. %
We substitute $\theta_{1}(x_{i}) = \phi - \alpha_{1}(x_{i})$ into \eqref{eqn:S-factor} and solve %
\begin{equation}
\label{eqn:step2}
  \min_{(\phi,\beta)} \sum_{i=1}^{M} \biggl| - \frac{S(\theta_{1}(x_{i});\phi,\beta)}{16\pi \beta \sqrt{r\,r_{1}(x_{i})}} - U_{i} \biggr|^{2}
\end{equation}
over the admissible set of $(\phi,\beta)$. %
Theorem \ref{thm:uniqueness} guarantees uniqueness of the minimizer in the noise-free case. %


\subsection{Recovery procedure for polygonal obstacles}
\label{subsec:polygon-recovery}

For a compact polygonal obstacle the local analysis near each vertex is unchanged, but the recorded field on $L$ is a superposition of geometric waves, primary diffracted waves from every visible vertex, and weaker multiple diffracted waves. %
As discussed in Section \ref{sec:introduction}, only visible vertices produce primary diffracted waves that reach $L$ without prior interaction with other parts of $\partial K$, whereas invisible vertices and vertices hit only by grazing rays are excluded. %

\begin{figure}[htb]
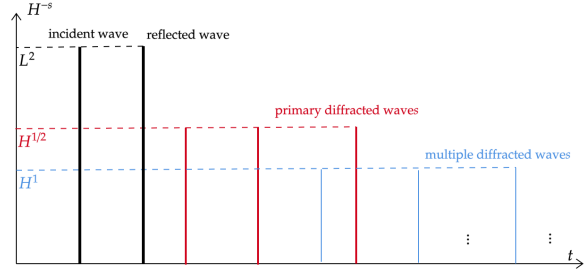

\centering
\tikz[>=latex,scale=1.0]{
\node at (6.4,0.6) {$\vdots$};
\node at (7.6,0.6) {$\vdots$};
\draw[thick,->] (0,0) -- (8.0,0) node[pos=1.0,right] {$t$};
\draw[thick,->] (0,0) -- (0,3.5) node[pos=1.0,right] {\scriptsize$H^{-s}$};
\draw[ultra thick] (1,0) -- ++(0,2.8) node[pos=1.0,above] {\scriptsize incident wave};
\draw[ultra thick] (2,0) -- ++(0,2.8) node[pos=1.0,above right] {\scriptsize reflected wave};
\draw[thick,dashed] (0,2.8) -- ++(2,0) node[pos=0.0,below right] {\scriptsize$L^{2}$};
\draw[very thick,red] (2.5,0) -- ++(0,1.8);
\draw[very thick,red] (3.5,0) -- ++(0,1.8);
\draw[very thick,red] (5.0,0) -- ++(0,1.8) node[pos=1.0,above] {\scriptsize primary diffracted waves};
\draw[thick,dashed,red] (0,1.8) -- ++(5.0,0) node[pos=0.0,below right] {\scriptsize$H^{1/2}$};
\draw[thick,blue] (4.5,0) -- ++(0,1.2);
\draw[thick,blue] (5.8,0) -- ++(0,1.2);
\draw[thick,blue] (7.0,0) -- ++(0,1.2) node[pos=1.0,above] {\scriptsize multiple diffracted waves};
\draw[thick,dashed,blue] (0,1.2) -- ++(7.0,0) node[pos=0.0,below right] {\scriptsize$H^{1}$};
}
\caption{Sobolev regularities of reflected and diffracted waves (schematic).}
\label{F_T_1}
\end{figure}

The constructive recovery of Theorem \ref{thm:polygon} proceeds in four stages, which we implement numerically as follows. %
First, one restricts attention to a sub-arc of $L$ on which primary diffracted waves from distinct vertices are pairwise separated in time and arrive after all geometric waves. %
Second, multiple diffracted waves are discarded by comparing Sobolev regularities: geometric waves have the largest singularities, followed by primary diffracted waves and then multiple diffracted waves; see Figure \ref{F_T_1}. %
Third, each remaining primary diffracted wave is localised using its diffracted arrival time $t_{j}(x) := |O-P_{j}| + |x-P_{j}|$. %
Fourth, Steps 1--2 of Section \ref{subsec:two-step} recover $(P_{j},\phi_{j},\beta_{j})$ from the corresponding arrival times \eqref{eqn:arrival_time} and amplitudes \eqref{eqn:jump}. %
Assembling the recovered corners gives the visible vertices and the local directions of their adjacent edges. %

After subtracting the incident field and isolating each primary diffracted signal, we apply Algorithm \ref{alg:two-step} to the corresponding visible vertex. %


\section{Numerical studies}
\label{sec:numerical}

We first describe the forward solver, then examine the simulated diffracted signals and present the reconstruction experiments. %


\subsection{Numerical method for the forward problem}
\label{subsec:forward-method}

The wave equation $( \partial_{tt} - \Delta ) u = 0$ is rewritten as the first-order velocity-pressure system %
\begin{equation}
\label{eqn:vp-system}
  \left\{
  \begin{aligned}
    \partial_{t} \bm{v} &= -\nabla p, \\
    \partial_{t} p &= -\nabla \cdot \bm{v},
  \end{aligned}
  \right.
  \qquad
    p := -\partial_{t} u,
  \quad
    \bm{v} := \nabla u = (v_{1},v_{2}).
\end{equation}
We combine a staggered-grid leapfrog scheme \cite{Virieux1986} with a classical split-field PML \cite{Berenger1994, Collino2001}. %
The method is second-order accurate in space and time. %


\subsubsection{Staggered grid and time stepping}

Pressure $p$ is stored at cell centres $(i+\tfrac{1}{2},j+\tfrac{1}{2})$, %
the horizontal velocity $v_{1}$ is stored at vertical faces $(i,j+\tfrac{1}{2})$, %
and the vertical velocity $v_{2}$ is stored at horizontal faces $(i+\tfrac{1}{2},j)$. %
Temporally, $p$ lives at integer steps $t^{n}$ and $\bm{v}$ at half-integer steps $t^{n+\frac{1}{2}}$. %
Away from PML and obstacles the leapfrog update reads %
\begin{equation*}
  \left\{
  \begin{aligned}
    \bm{v}^{n+\frac{1}{2}}
      &= \bm{v}^{n-\frac{1}{2}} - \Delta t \,\nabla_{h} p^{n}, \\
    p^{n+1}
      &= p^{n} - \Delta t \,\nabla_{h} \cdot \bm{v}^{n+\frac{1}{2}}.
  \end{aligned}
  \right.
\end{equation*}
Here $\nabla_{h}$ denotes centred differences of spacing $h$. %
The scheme is conditionally stable under a standard CFL restriction \cite{Courant1928, Strikwerda2004}. %


\subsubsection{Initialisation}

Given $u|_{t=0} = \varphi$ and $\partial_{t} u|_{t=0} = \psi$, set %
\begin{equation*}
  p^{0} = -\psi,
  \qquad
  \bm{v}^{0} = \nabla_{h} \varphi,
\end{equation*}
and advance velocity to the first half step by a second-order Taylor step, %
\begin{equation*}
  \bm{v}^{\pm \frac{1}{2}}
    = \nabla_{h} \varphi \pm \tfrac{1}{2} \Delta t \,\nabla_{h} \psi.
\end{equation*}
In the pulse experiments below we take $\varphi = 0$ and approximate $\psi = \delta_{O}$ by a narrow Gaussian centred at the source $O$. %


\subsubsection{Standard split-field PML}

Inside the absorbing layers the fields are split direction-wise, $p = \sum_{k\in\{1,2\}} p_{k}$, with %
\begin{equation*}
  \left\{
  \begin{aligned}
    \partial_{t} v_{k} + \sigma_{k} v_{k} &= -\partial_{k} p, \\
    \partial_{t} p_{k} + \sigma_{k} p_{k} &= -\partial_{k} v_{k},
  \end{aligned}
  \right.
\end{equation*}
and polynomial damping $\sigma(s) = \sigma_{\max}(s/L_{\mathrm{PML}})^{m}$. %
A semi-implicit discretisation yields %
\begin{equation*}
  \left\{
  \begin{aligned}
    v_{k}^{n+\frac{1}{2}}
      &= e^{-\sigma_{k}\Delta t} \, v_{k}^{n-\frac{1}{2}}
       - \sigma_{k}^{-1} \bigl( 1 - e^{-\sigma_{k}\Delta t} \bigr) \, \partial_{k,h} p^{n}, \\
    p_{k}^{n+1}
      &= e^{-\sigma_{k}\Delta t} \, p_{k}^{n}
       - \sigma_{k}^{-1} \bigl( 1 - e^{-\sigma_{k}\Delta t} \bigr) \, \partial_{k,h} v_{k}^{n+\frac{1}{2}}.
  \end{aligned}
  \right.
\end{equation*}
When $\sigma_{k} = 0$ one recovers the interior leapfrog update. %
The physical pressure is $p^{n} = \sum_{k\in\{1,2\}} p_{k}^{n}$. %


\subsubsection{Recovery of the wave field $u$}

If the scalar field $u$ is required, integrate $\partial_{t} u = -p$ by the trapezoidal rule, %
\begin{equation}
  u^{n+1} = u^{n} - \tfrac{1}{2} \Delta t \, \bigl( p^{n} + p^{n+1} \bigr),
  \qquad u^{0} = \varphi.
\end{equation}


\subsubsection{Neumann boundary conditions via $v_{1}$ and $v_{2}$}
\label{subsec:neumann-vxvy}

On a rigid obstacle the Neumann condition $\partial_{\nu} u = 0$ is equivalent to vanishing normal velocity, %
\begin{equation}
\label{eqn:neumann-v}
  \nu \cdot \bm{v} = \nu_{1} v_{1} + \nu_{2} v_{2} = 0 \qquad \text{on } \partial K,
\end{equation}
where $\nu = (\nu_{1},\nu_{2})$ is the unit outward normal. 

\paragraph{Axis-aligned edges}
If an edge of $\partial K$ is vertical, then $\nu = (\pm 1,0)$ and \eqref{eqn:neumann-v} reduces to $v_{1} = 0$. %
If the edge is horizontal, then $\nu = (0,\pm 1)$ and $v_{2} = 0$. %
After each velocity update, the corresponding face values of $v_{1}$ or $v_{2}$ that lie on such edges are overwritten by zero. %
Cells marked as interior to $K$ keep $\bm{v} \equiv 0$ and do not update $p$. %

\paragraph{Slanted polygonal edges}
For a general edge with outward unit normal $\nu$, we consider the set of staggered velocity nodes whose control faces intersect that edge. %
At each such node, one retains only the tangential component of velocity: %
\begin{equation}
\label{eqn:project-v}
\bm{v} ~\longleftarrow~ \bm{v} - ( \nu \cdot \bm{v} ) \nu.
\end{equation}
When the edge is axis-aligned, this projection collapses to setting the single normal component to zero, as above. %
Ghost values of $p$ inside $K$, when needed for the centred difference of $\nabla_{h} p$ next to the boundary, are obtained by even reflection across the interface so that the discrete normal derivative of $p$ remains consistent with \eqref{eqn:neumann-v}; cf. \cite{Gibou2002, Coco2013, Arias2018}. %

The combination of (i) solid-cell masking, (ii) zeroing or projecting $(v_{1},v_{2})$ on boundary faces, and (iii) reflected ghost pressures yields a stable second-order treatment of Neumann conditions for polygonal obstacles on the staggered grid used in all experiments below. %


\subsection{Numerical simulation of the diffracted wave}
\label{subsec:diff-sim}

Before addressing inversion, we validate that the forward solver reproduces the qualitative structure predicted by the Cheeger--Taylor theory recalled in Section \ref{sec:reconstruction-scheme}: a geometric wavefront $\mathcal{G}_{t}^{O}$, a diffracted wavefront $\mathcal{D}_{t}^{O}$ arriving at $t_{P}(x) = r + r_{1}(x)$, and angular sectors in which reflection may be present or absent while diffraction persists. %
We take a single convex corner as the model problem. %
Throughout this subsection a full-circle receiver configuration is used only to visualise wave propagation and signal structure. %
The source is approximated by a Gaussian pulse of width $\sigma = 0.02$, %
\begin{equation*}
  \mathcal{S}_{\sigma}(x;O)
    = \frac{1}{2\pi \sigma^{2}} \exp\biggl( -\frac{\lvert x-O \rvert^{2}}{2 \sigma^{2}} \biggr).
\end{equation*}
The computational window is $[-1,1]^{2}$ with PML layers, mesh size $h = 2\times10^{-3}$, and a CFL-compliant time step. %
The opening parameter is $\beta = \tfrac{5}{6} \in (\tfrac{1}{2},1)$. %


\subsubsection{Convex corner: configuration and snapshots}

\begin{figure}[htb]
\centering
\includegraphics[width=0.28\textwidth,trim=40 0 70 0,clip]{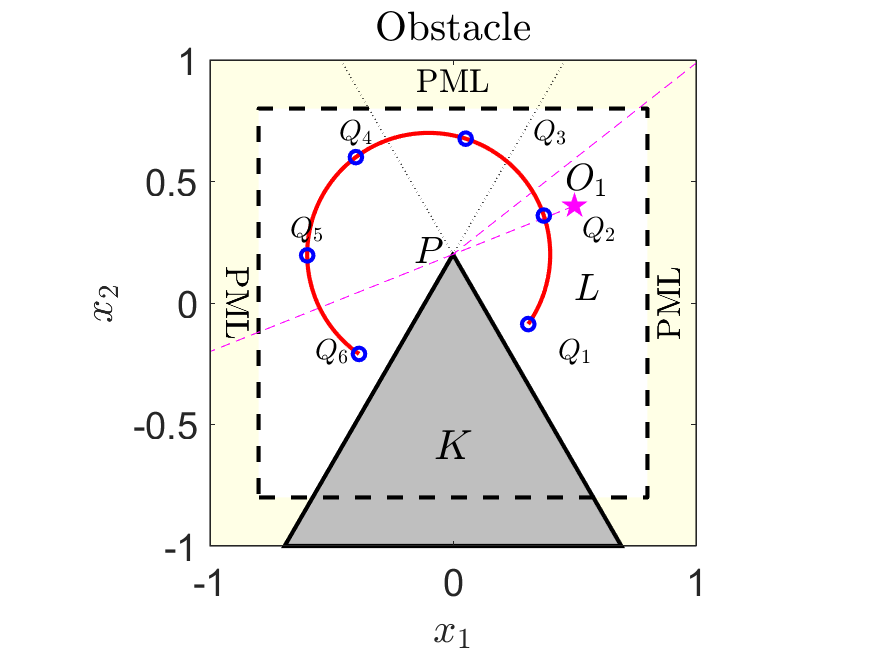}
\caption{Convex corner with source $O_{1}$ (full-circle receivers). %
Black dashed lines: PML; black solid lines: obstacle; magenta star: source; red arcs: receivers.}
\label{fig:obstacle_O1}
\end{figure}

\begin{figure}[htb]
\centering
\setlength{\tabcolsep}{4pt}
\begin{tabular}{cccc}
\includegraphics[width=0.23\textwidth,trim=10 0 30 0,clip]{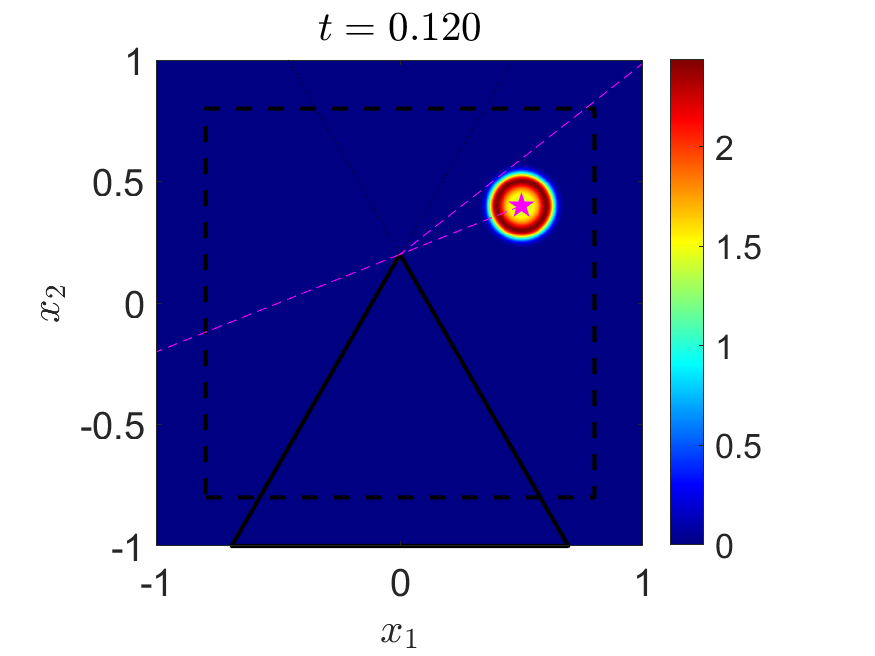} &
\includegraphics[width=0.23\textwidth,trim=10 0 30 0,clip]{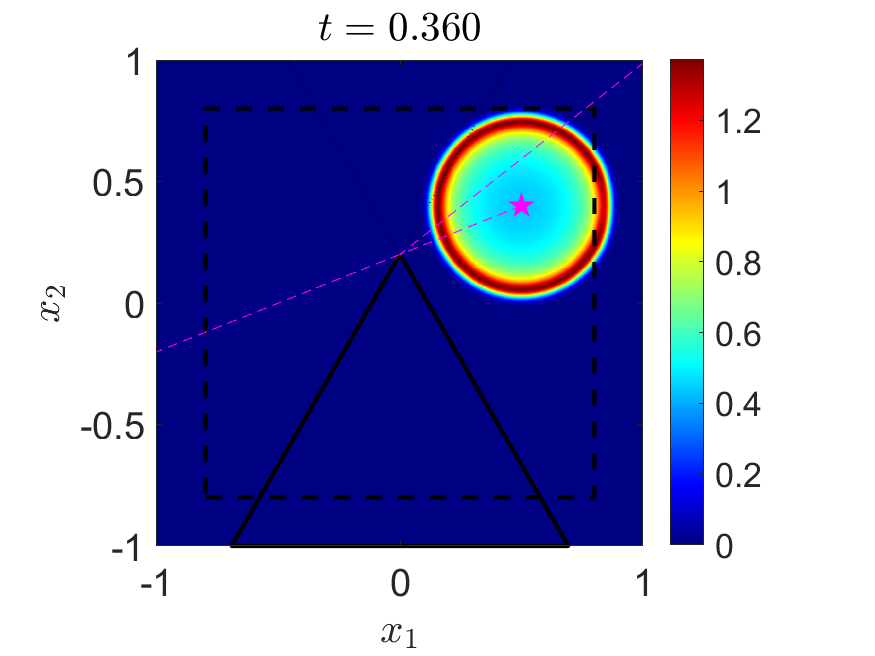} &
\includegraphics[width=0.23\textwidth,trim=10 0 30 0,clip]{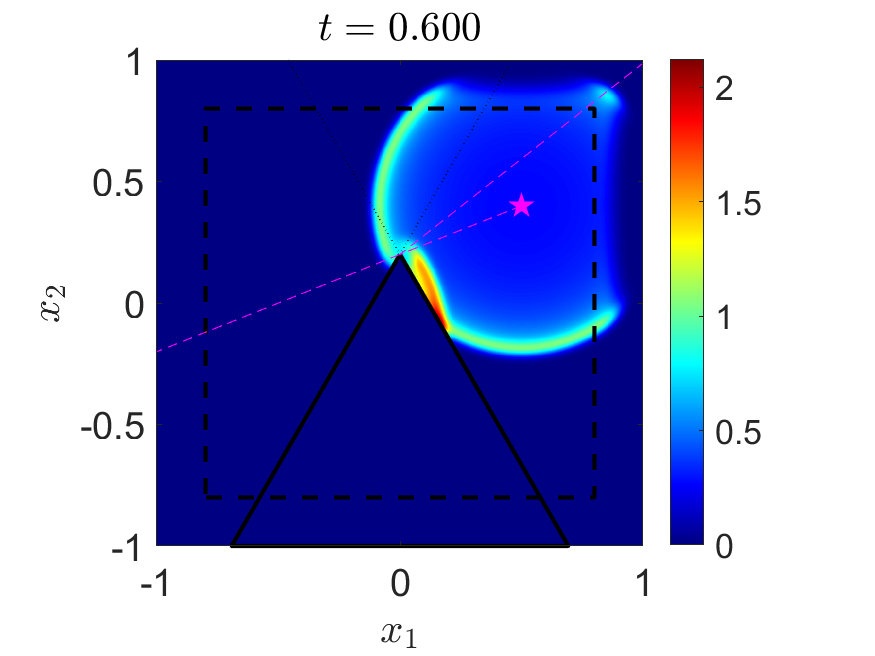} &
\includegraphics[width=0.23\textwidth,trim=10 0 30 0,clip]{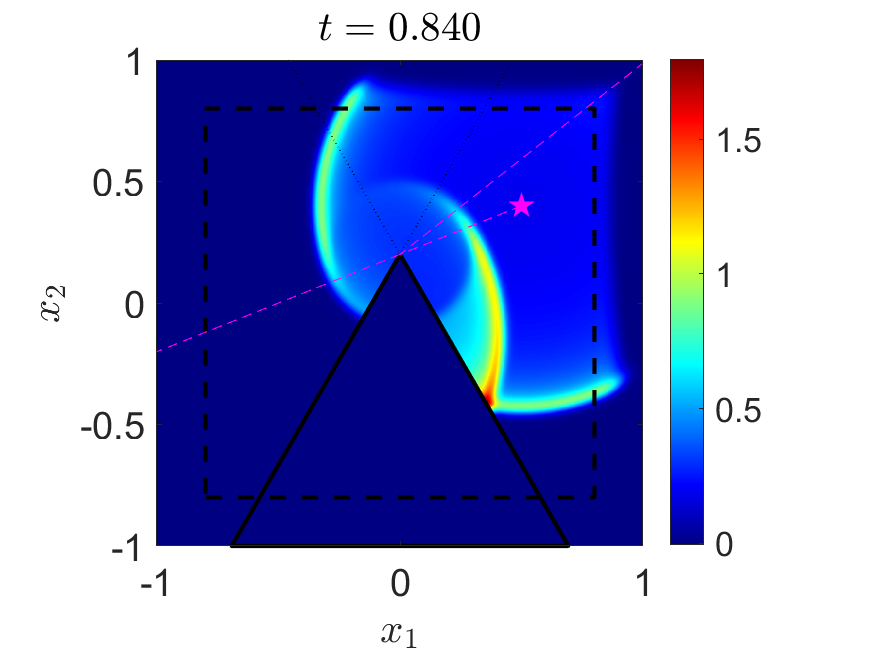} \\[-0.4ex]
{\scriptsize (1)} & {\scriptsize (2)} & {\scriptsize (3)} & {\scriptsize (4)} \\[0.6ex]
\includegraphics[width=0.23\textwidth,trim=10 0 30 0,clip]{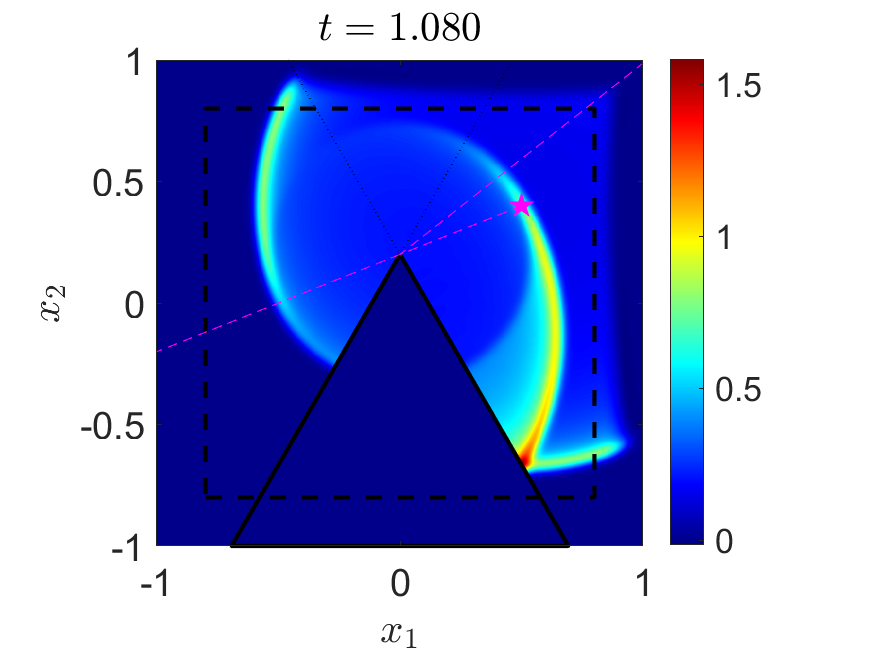} &
\includegraphics[width=0.23\textwidth,trim=10 0 30 0,clip]{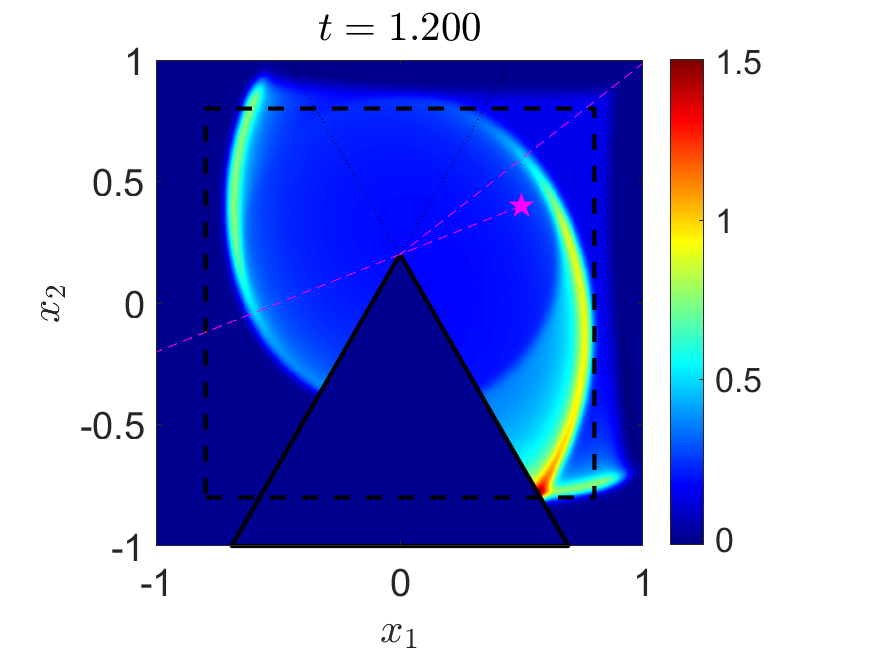} &
\includegraphics[width=0.23\textwidth,trim=10 0 30 0,clip]{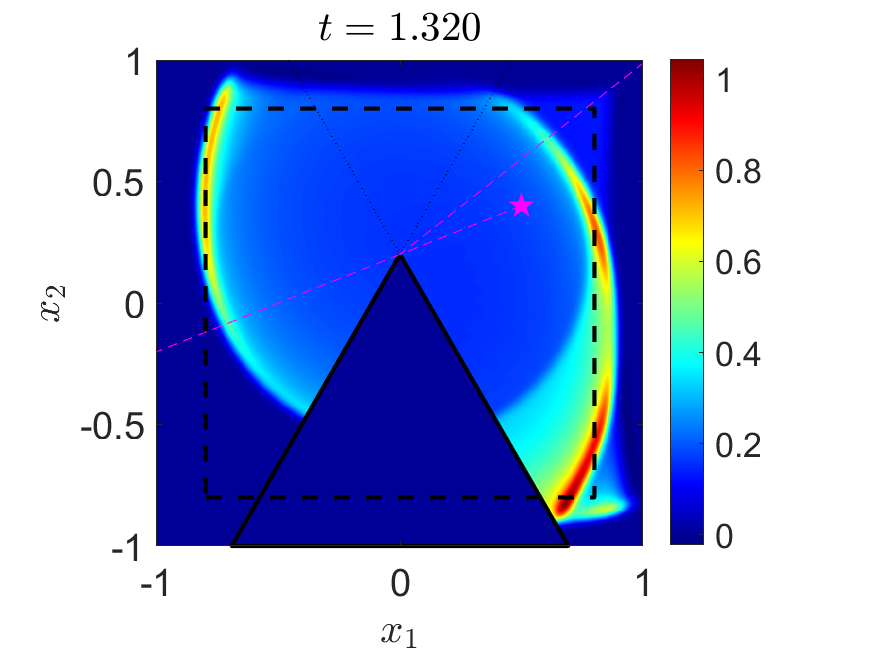} &
\includegraphics[width=0.23\textwidth,trim=10 0 30 0,clip]{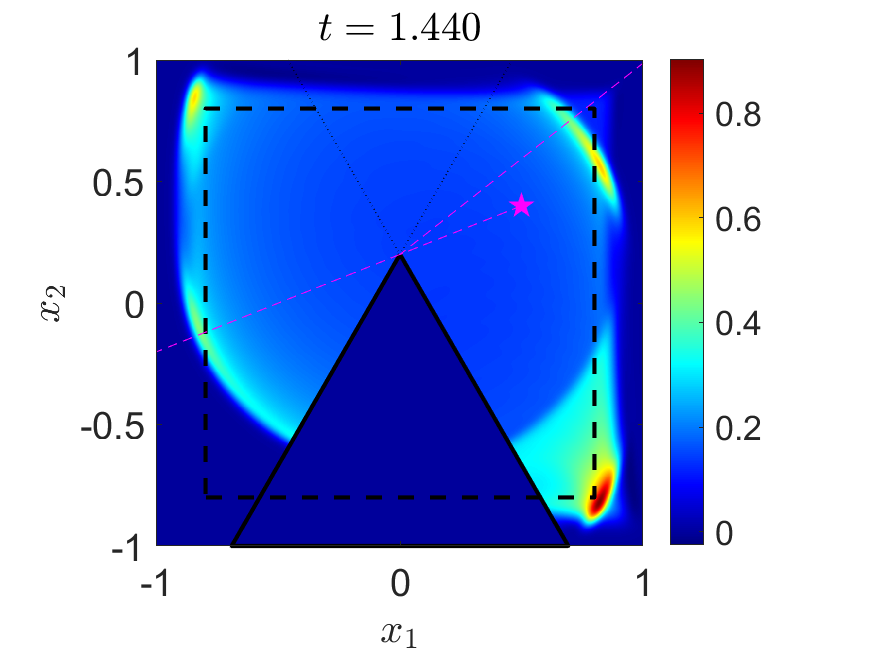} \\[-0.4ex]
{\scriptsize (5)} & {\scriptsize (6)} & {\scriptsize (7)} & {\scriptsize (8)} \\
\end{tabular}
\caption{Snapshots of the total wave field for the convex corner ($\sigma = 0.02$). %
Reflected wavefronts on the faces and a circular diffracted wavefront from $P$ are clearly visible.}
\label{fig:obstacle_O1_u}
\end{figure}

Figure \ref{fig:obstacle_O1} shows the simulation geometry: a single convex corner with tip $P$, a point source at $O_{1}$, and a receiver curve in the exterior domain (red arcs) surrounding the tip. %
Figure \ref{fig:obstacle_O1_u} displays successive snapshots of the total wave field $u$. %
One observes, in chronological order, %
(i) the free expansion of the incident pulse, %
(ii) specular reflection on the illuminated faces, %
and (iii) a circular diffracted wavefront emanating from $P$ after the incident wave has reached the tip, %
in agreement with the Cheeger--Taylor / GTD description recalled in Section \ref{sec:reconstruction-scheme}. %

\subsubsection{Regions I, II, III and angular Regions A, B, C}

\begin{figure}[htb]
\centering
\begin{tabular}{cc}
\resizebox{0.36\textwidth}{!}{
\tikz[>=latex]{
\node at (0,0) {\includegraphics[width=0.5\textwidth,trim=10 0 30 0,clip]{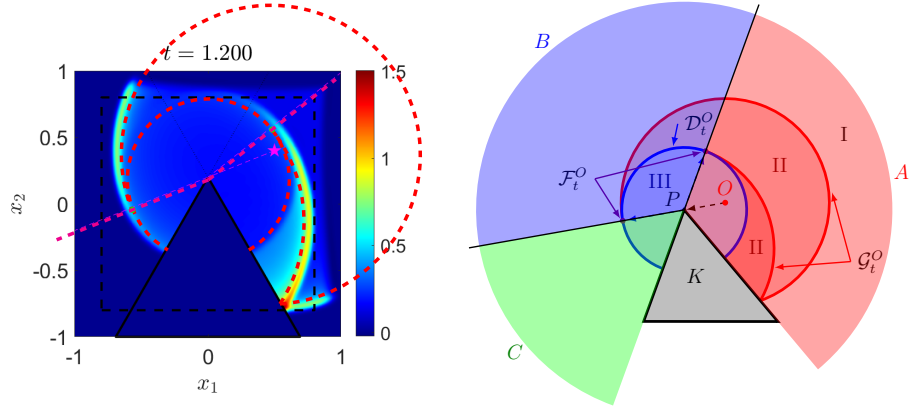}};
\coordinate (P) at (-0.14,0.61);
\coordinate (O) at (0.87,1.02);
\coordinate (OO) at (-1,-0.07);
\coordinate (A) at ($ (P) + (-60:1.3) $);
\draw[ultra thick,red,dashed] (A) arc (-60:240:1.3);
\coordinate (B) at ($ (P) + (-60:2.3) $);
\draw[ultra thick,red,dashed] (B) arc (-86.5:201.8:2.4);
\draw[ultra thick,red,dashed] (B) arc (-33:38.3:2.4);
\coordinate (C) at ($ (O) + (201.8:2.4) $);
\draw[ultra thick,green] (C) circle (0.02);
\coordinate (D) at ($ (OO) + (38.3:2.4) $);
\draw[ultra thick,green] (D) circle (0.02);
\draw[ultra thick,magenta,dashed] (P) -- ($ (C) !-1! (P) $);
\draw[ultra thick,magenta,dashed] (P) -- ($ (D) !-1! (P) $);
}}
&
\resizebox{0.36\textwidth}{!}{
\tikz[>=latex,rotate=10,scale=1]{
\coordinate (O)  at (0,0);
\coordinate (P)  at (-1,0);
\coordinate (R)  at ($ (P) + (-120:1) $);
\coordinate (F1) at (-2.5,0);
\coordinate (F2) at ($ (P) + (60:1.5) $);
\draw[ultra thick,fill=red!20!white,draw=red] (O) -- ($ (O) + (-120:2.5) $) arc (-120:180:2.5) -- cycle;
\draw[ultra thick,fill=red!40!white,draw=red] (R) -- ($ (R) + (-45:2.5) $) arc (-45:60:2.5) -- cycle;
\draw[ultra thick,fill=blue!20!white,draw=blue] (P) circle (1.5);
\draw[ultra thick,fill=lightgray] (P) -- ($ (P) + (-60:3.5) $) -- ($ (P) + (-120:2.85) $) -- cycle;
\draw[thick,draw=black,dashed,->] (O) -- ($ (P)  !0.05! (O) $);
\draw[thick,draw=blue,dashed,->] (P) -- ($ (F1) !0.07! (P) $);
\draw[thick,draw=blue,dashed,->] (P) -- ($ (F2) !0.07! (P) $);
\node (G) at (3.2,-2) {$\mathcal{G}_{t}^{O}$};
\draw[thick,draw=red,->] (G.west) -- ($ (2.5,0) !0.1! (G) $);
\draw[thick,draw=red,->] (G.west) -- ($ (R) + (-20:2.55) $);
\node (D) at ($ (P) + (70:2.1) $) {$\mathcal{D}_{t}^{O}$};
\draw[thick,draw=blue,->] (D.west) -- ($ (P) + (90:1.55) $);
\node (F) at (-3.5,1.2) {$\mathcal{F}_{t}^{O}$};
\draw[thick,draw=red!50!blue,->] (F.east) -- ($ (F1) !0.05! (F) $);
\draw[thick,draw=red!50!blue,->] (F.east) -- ($ (F2) !0.03! (F) $);
\fill[red] (O) node[above] {$O$} circle (0.07);
\fill[black] (P) node[above left=-2] {$P$} circle (0.05);
\fill[red!50!blue] (F1) circle (0.07);
\fill[red!50!blue] (F2) circle (0.07);
\node at ($ (P) + (-90:1.6) $) {$K$};
\node at ($ (O) + ( 20:3.3) $) {\textrm{I}};
\node[red!50!black] at ($ (O) + ( 25:1.6) $) {\textrm{II}};
\node[red!50!black] at ($ (O) + (-65:1.3) $) {\textrm{II}};
\node[blue!50!black] at ($ (P) + (120:0.9) $) {\textrm{III}};
\fill[red,opacity=0.35] (P) -- ($ (P) + (-60:5) $) arc (-60:60:5) -- cycle;
\fill[blue,opacity=0.35] (P) -- ($ (P) + ( 60:5) $) arc (60:180:5) -- cycle;
\fill[green,opacity=0.35] (P) -- ($ (P) +(-120:5) $) arc (-120:-180:5) -- cycle;
\draw[thick] (P) -- ($ (F1) !-2.5! (P) $);
\draw[thick] (P) -- ($ (F2) !-2.5! (P) $);
\node[red] at ($ (P) + (0:5.3) $) {$A$};
\node[blue] at ($ (P) + (120:5.3) $) {$B$};
\node[green!50!black] at ($ (P) + (210:5.3) $) {$C$};
}}
\end{tabular}
\caption{Left: late-time snapshot with geometric and diffracted wavefronts (red) and forward rays from $P$ (magenta). %
Right: time Regions I, II, III and angular Regions A, B, C for diffraction by a convex corner.}
\label{fig:regions-ABC}
\end{figure}

Recall the causal partition of Section \ref{sec:reconstruction-scheme}: no wave has reached Region I, Region II has been swept by geometric waves only, and Region III has been swept by the diffracted wavefront $t_{P}(x) = r + r_{1}(x)$. %
Independent of this time partition, the exterior near $P$ is divided into three angular sectors A, B, and C by the two lines extended beyond $P$ and the forward set $\mathcal{F}_{t}^{O}$; see Figure \ref{fig:regions-ABC}. %
The left panel of Figure \ref{fig:regions-ABC} overlays geometric and diffracted wavefronts on a late-time snapshot: the red dashed curves mark those wavefronts, and the magenta dashed lines mark the forward rays separating the angular sectors. %
The right panel sketches the corresponding time Regions I, II, III and angular Regions A, B, C. %
Region A is the illuminated sector in which receivers typically record a strong pulse and a specularly reflected signal, together with a weaker diffracted signal. %
Region B is an intermediate sector in which the reflected signal is absent, so after the incident signal, the diffracted signal dominates. %
Region C is the sector that the incident wave does not reach, where only diffracted waves are present. %
This angular classification, consistent with the geometric theory, explains why diffracted signals alone remain informative even when a receiver fails to record reflected signals. %

\subsubsection{Full-circle space-time records and pointwise signals}

\begin{figure}[htb]
\centering
\includegraphics[width=0.3\textwidth]{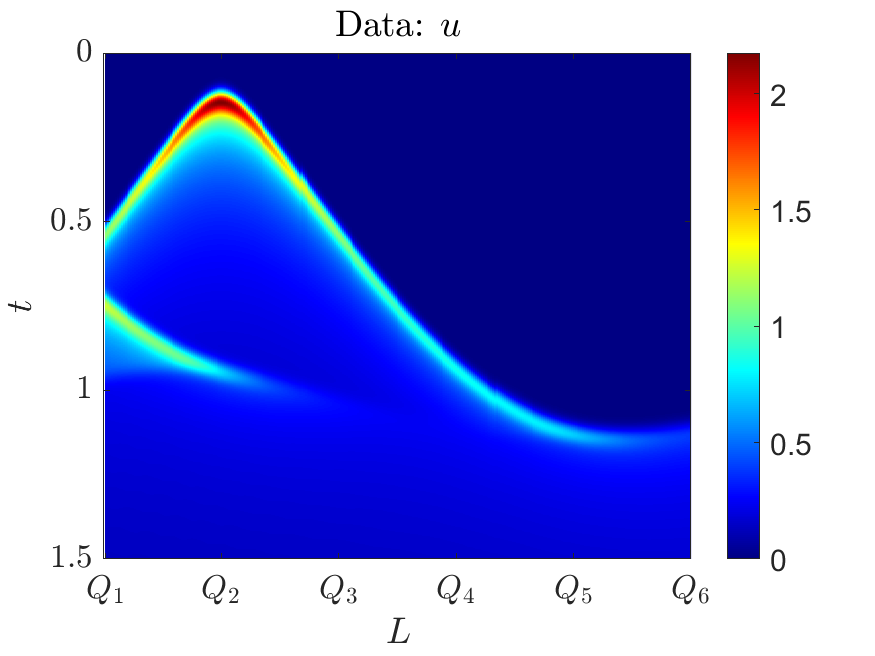}
\includegraphics[width=0.3\textwidth]{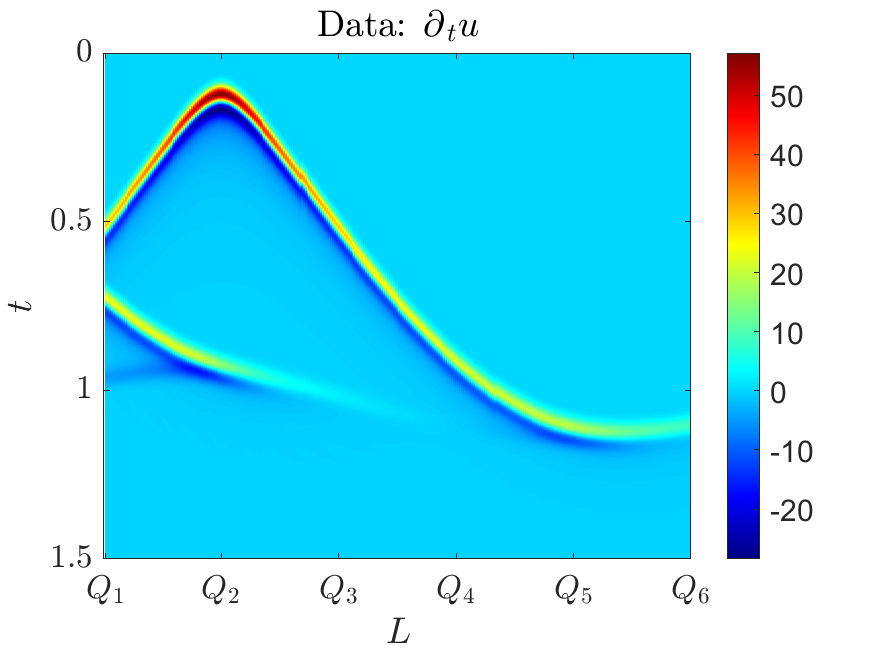}\\[0.5ex]
\includegraphics[width=0.3\textwidth]{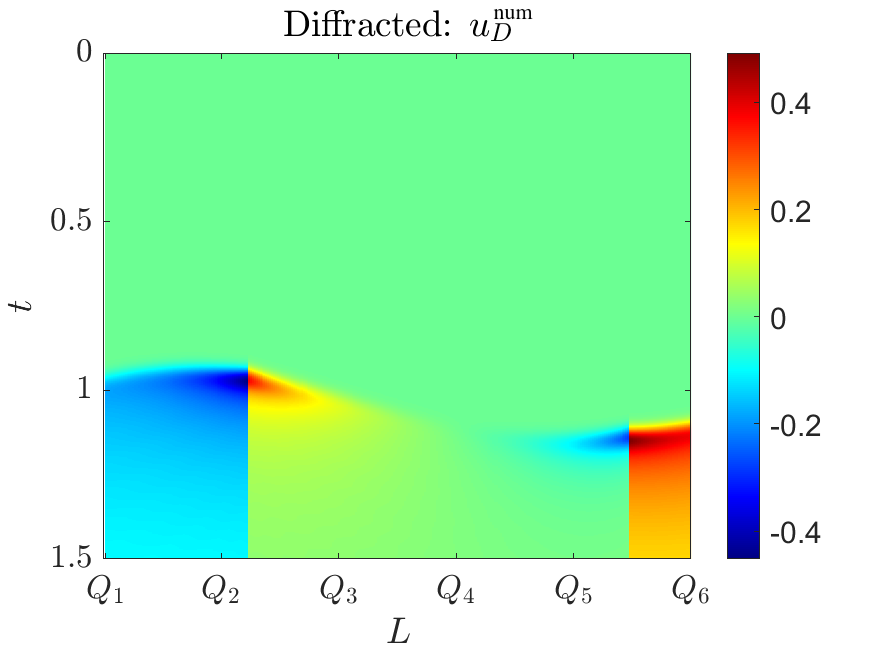}
\includegraphics[width=0.3\textwidth]{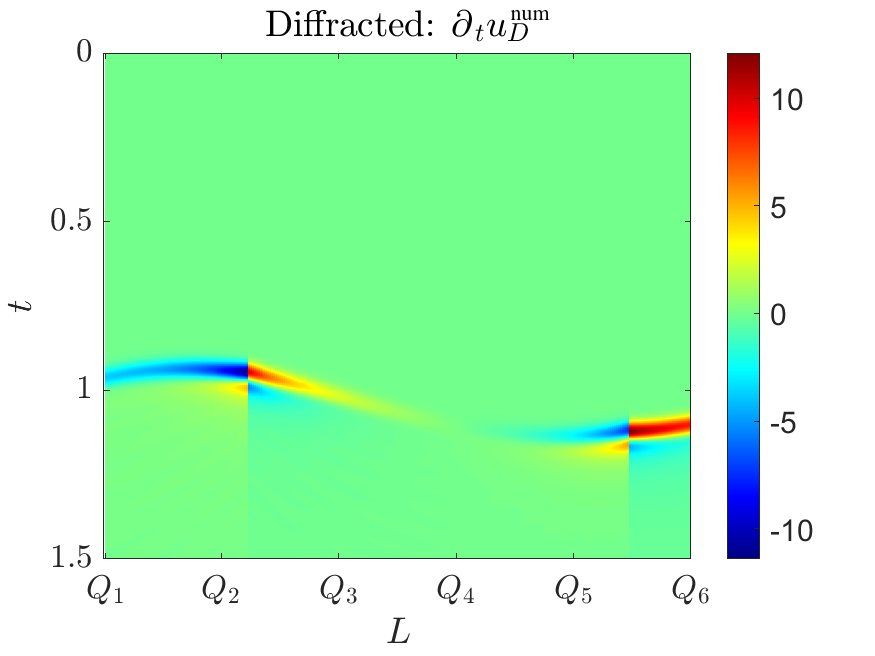}\\
\caption{Full-circle records: total wave field $u$ and $\partial_{t} u$ (top); %
extracted diffracted signal $u_{\mathcal{D}}^{\mathrm{num}}$ and its time derivative (bottom). %
Horizontal axis: receiver angle; vertical axis: time.}
\label{fig:data-full}
\end{figure}

\begin{figure}[htb]
\centering
\includegraphics[width=0.32\textwidth,trim=10 0 30 0,clip]{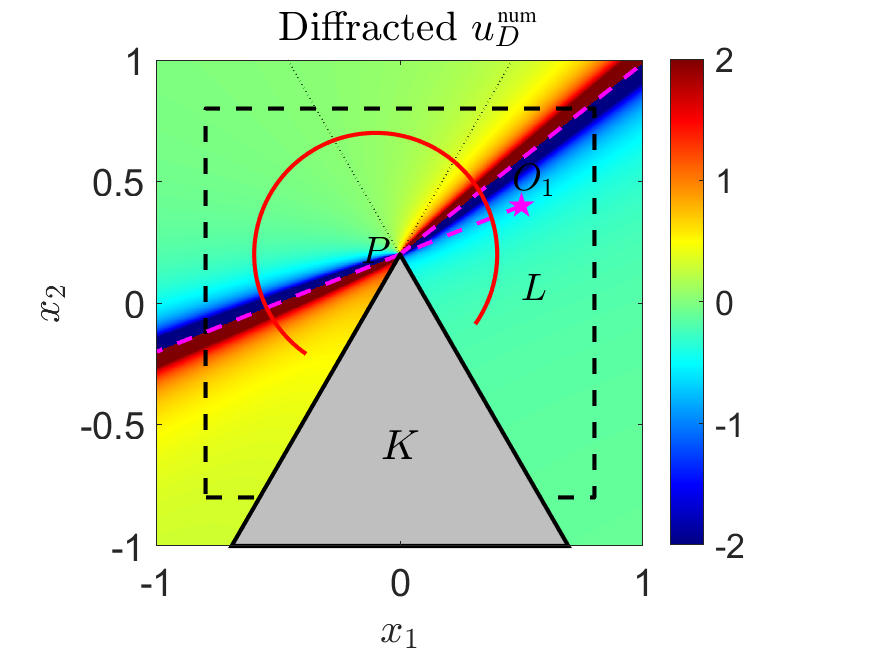}
\caption{Amplitude distribution of the extracted diffracted signal on the full circle (forward simulation only).}
\label{fig:uD-amp-portrait}
\end{figure}

Figure \ref{fig:data-full} shows space-time plots on the full circle of receivers: the total wave field $u$, its time derivative, and the extracted diffracted signal $u_{\mathcal{D}}^{\mathrm{num}}$. %
The incident and reflected signals appear as bands of large amplitude in the space-time plots. %
Where geometric signals are present, the diffracted signal arrives later, near $t_{P}(x) = r + r_{1}(x)$. %
Figure \ref{fig:uD-amp-portrait} displays the amplitude distribution of the extracted diffracted signal, which already exhibits the angular modulation predicted by $S(\theta_{1}(x);\phi,\beta)$ in \eqref{eqn:diffract}. %

\begin{figure}[htb]
\centering
\includegraphics[width=0.32\textwidth,trim=0 20 0 65,clip]{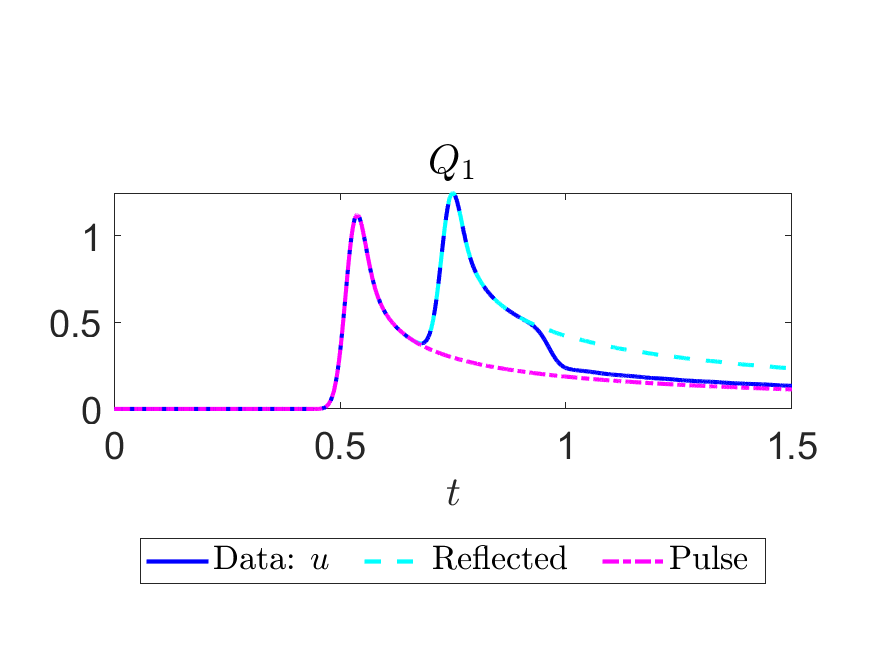}
\includegraphics[width=0.32\textwidth,trim=0 20 0 65,clip]{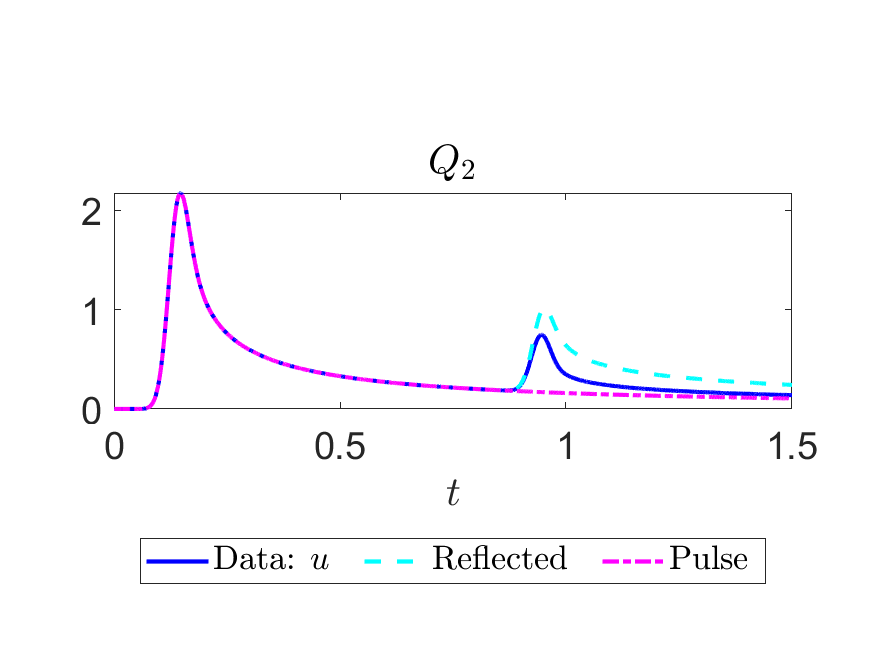}\\
\includegraphics[width=0.32\textwidth,trim=0 20 0 65,clip]{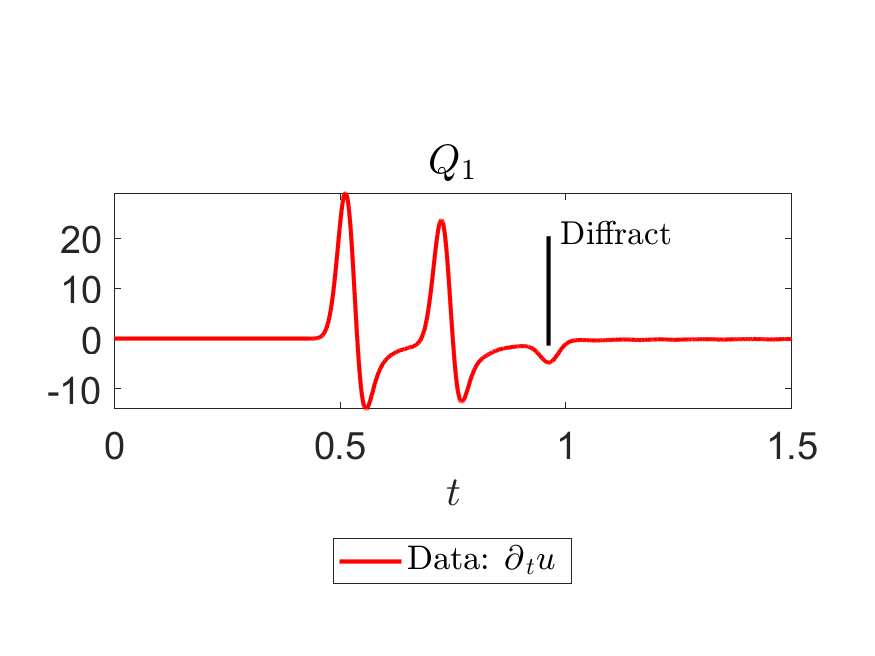}
\includegraphics[width=0.32\textwidth,trim=0 20 0 65,clip]{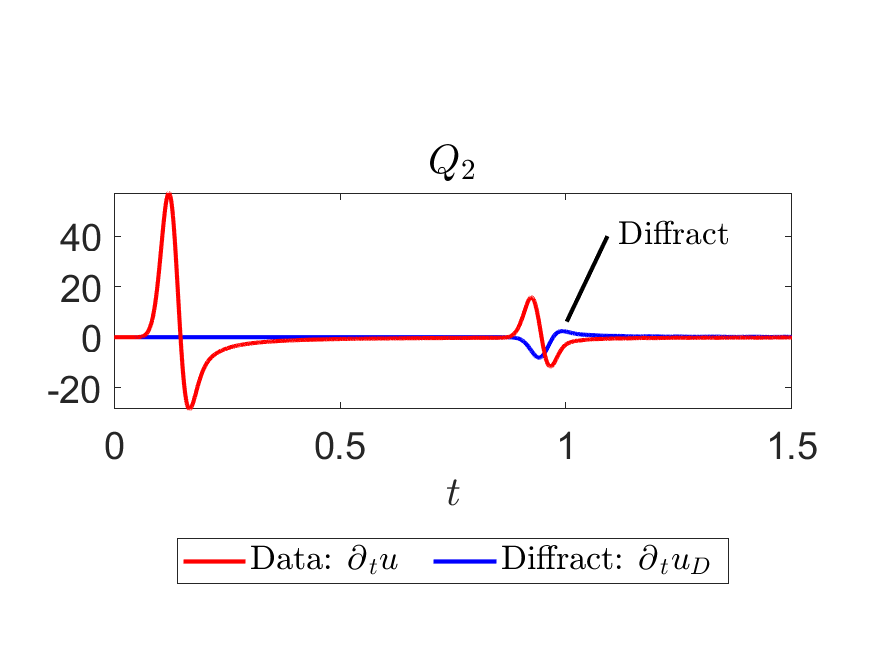}\\
\includegraphics[width=0.32\textwidth,trim=0 20 0 65,clip]{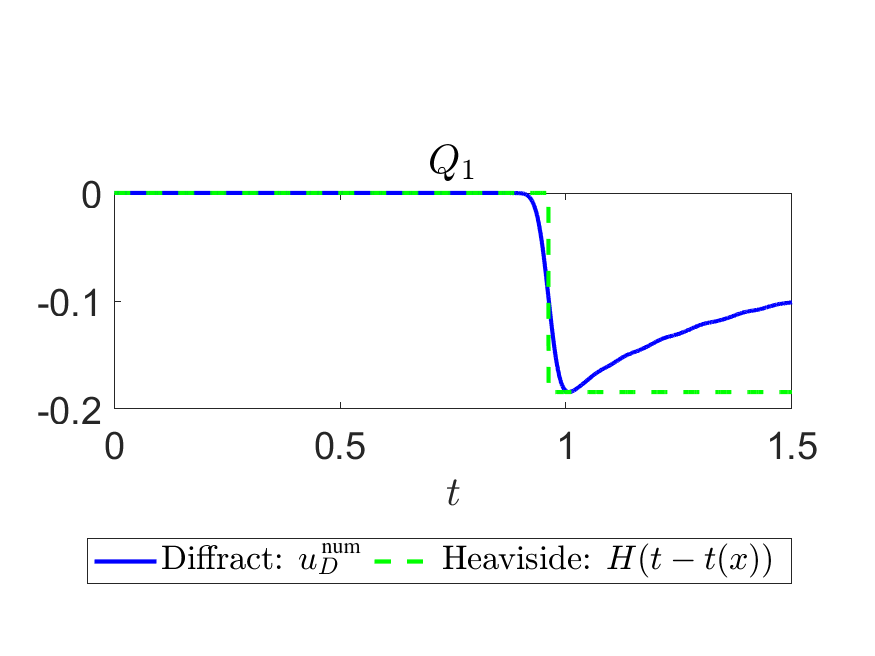}
\includegraphics[width=0.32\textwidth,trim=0 20 0 65,clip]{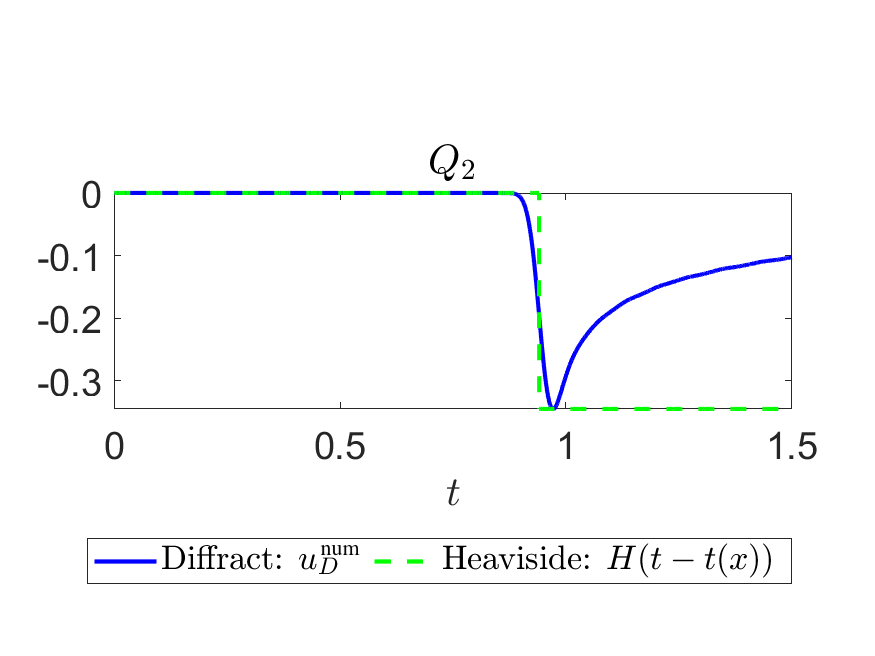}\\
\caption{Region A: pointwise signals at positions $Q_{1}$ and $Q_{2}$. %
Rows: $u$, $\partial_{t} u$, and extracted $u_{\mathcal{D}}^{\mathrm{num}}$. %
Three signals (incident, reflected, diffracted) are visible.}
\label{fig:signal-A}
\end{figure}

\begin{figure}[htb]
\centering
\includegraphics[width=0.3\textwidth,trim=0 20 0 65,clip]{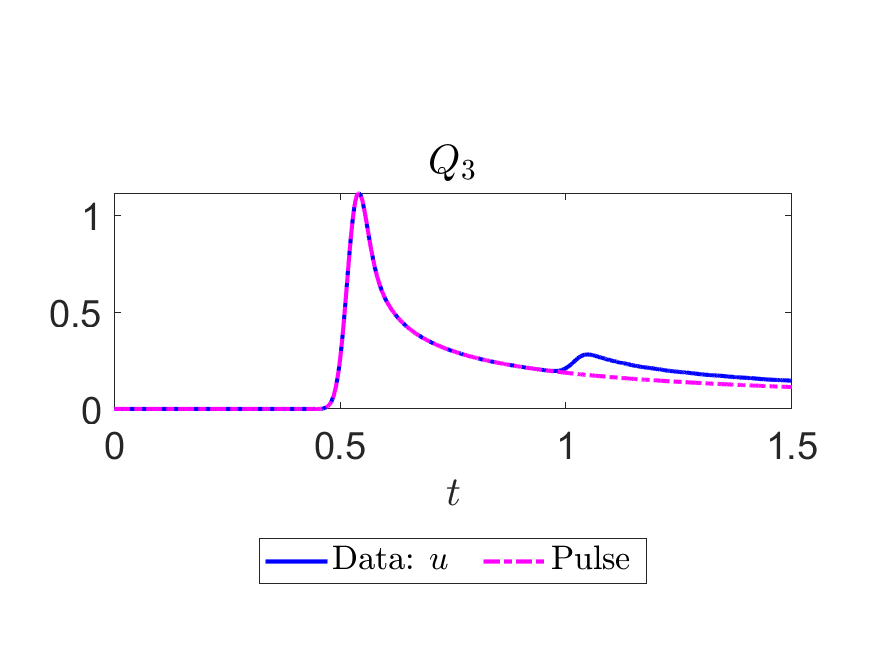}
\includegraphics[width=0.3\textwidth,trim=0 20 0 65,clip]{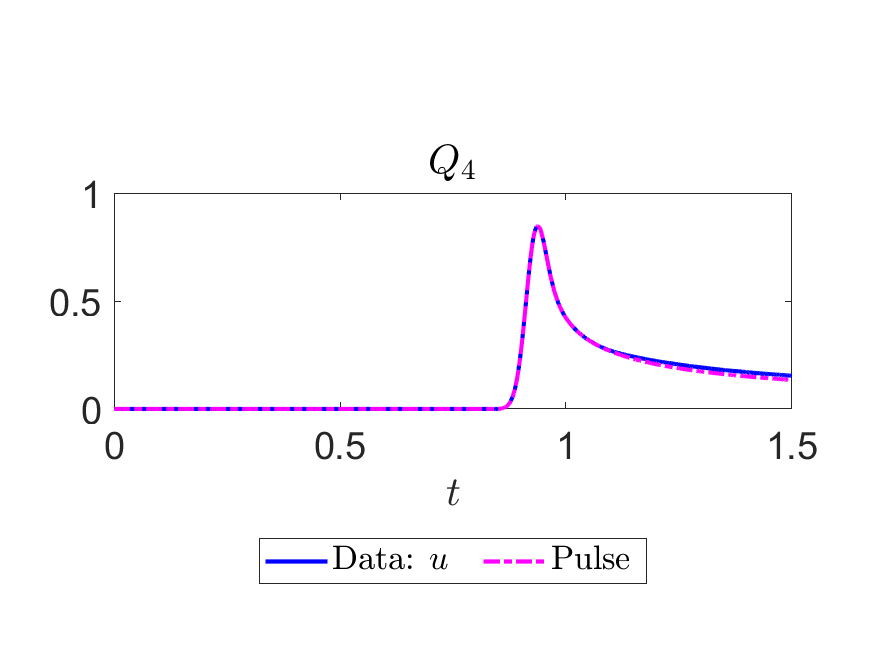}
\includegraphics[width=0.3\textwidth,trim=0 20 0 65,clip]{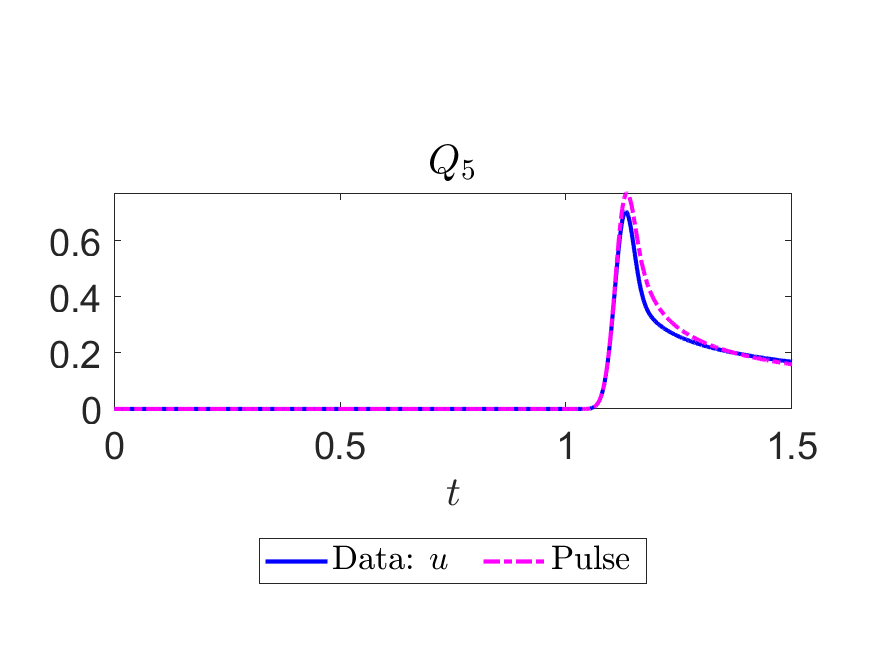}\\
\includegraphics[width=0.3\textwidth,trim=0 20 0 65,clip]{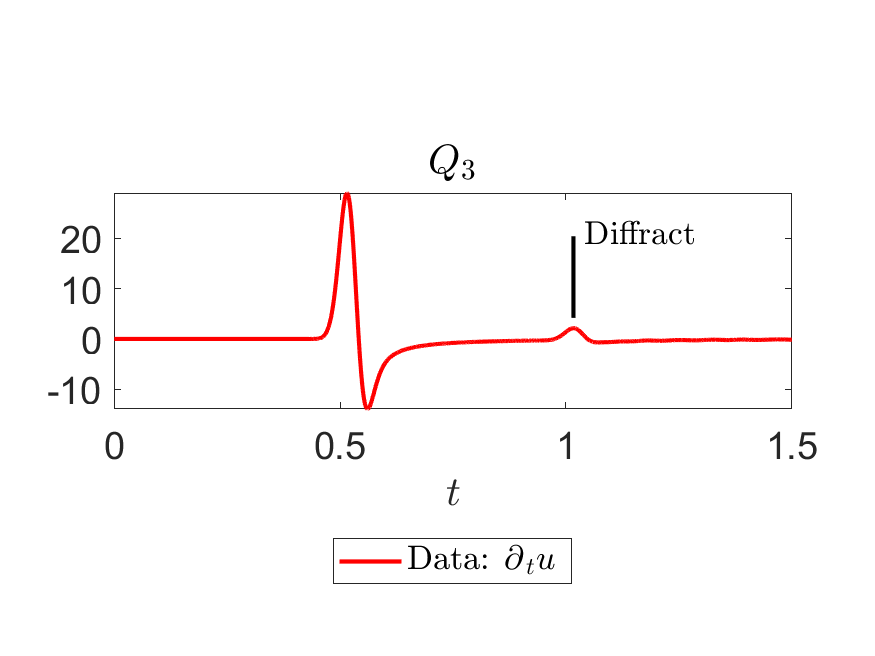}
\includegraphics[width=0.3\textwidth,trim=0 20 0 65,clip]{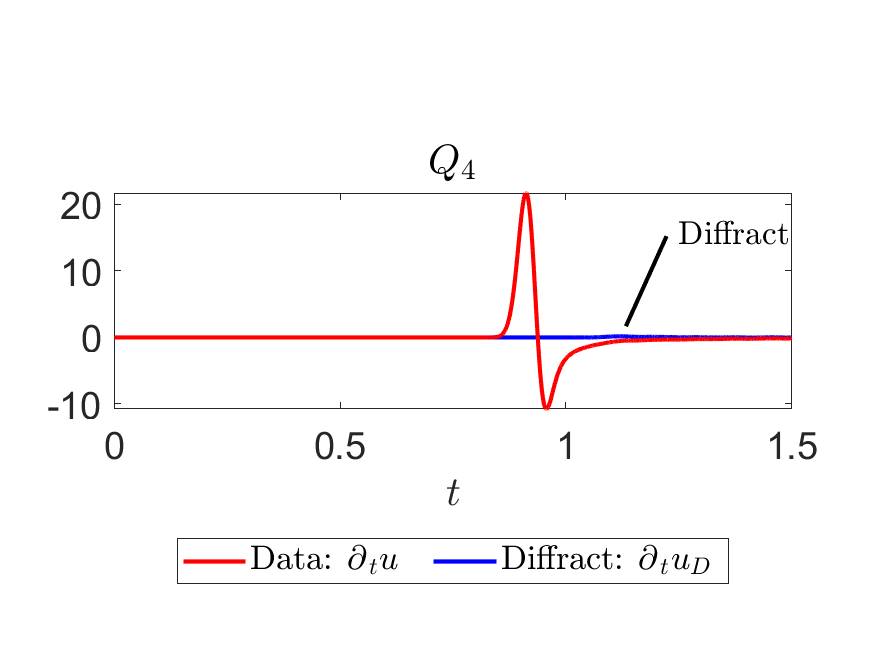}
\includegraphics[width=0.3\textwidth,trim=0 20 0 65,clip]{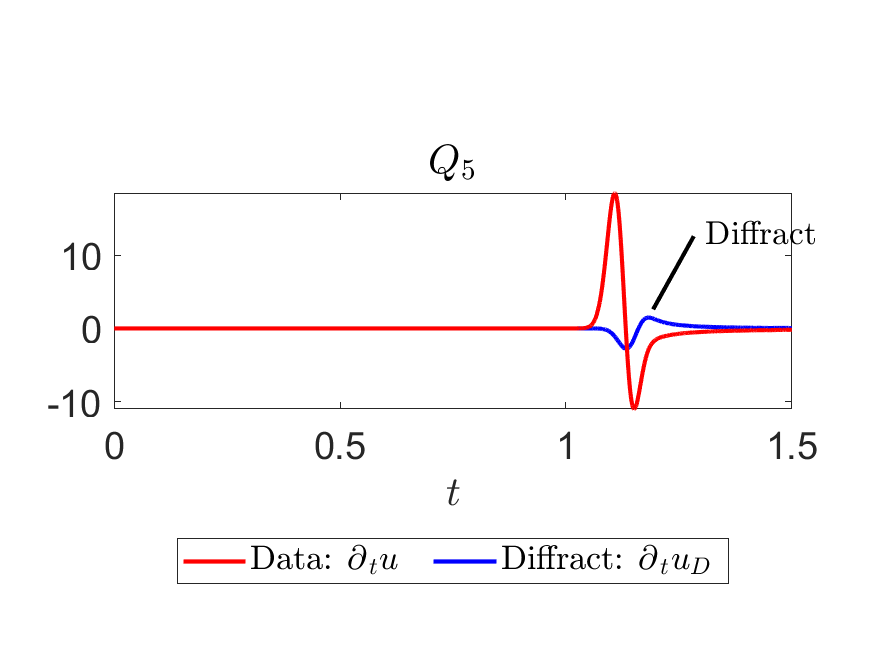}\\
\includegraphics[width=0.3\textwidth,trim=0 20 0 65,clip]{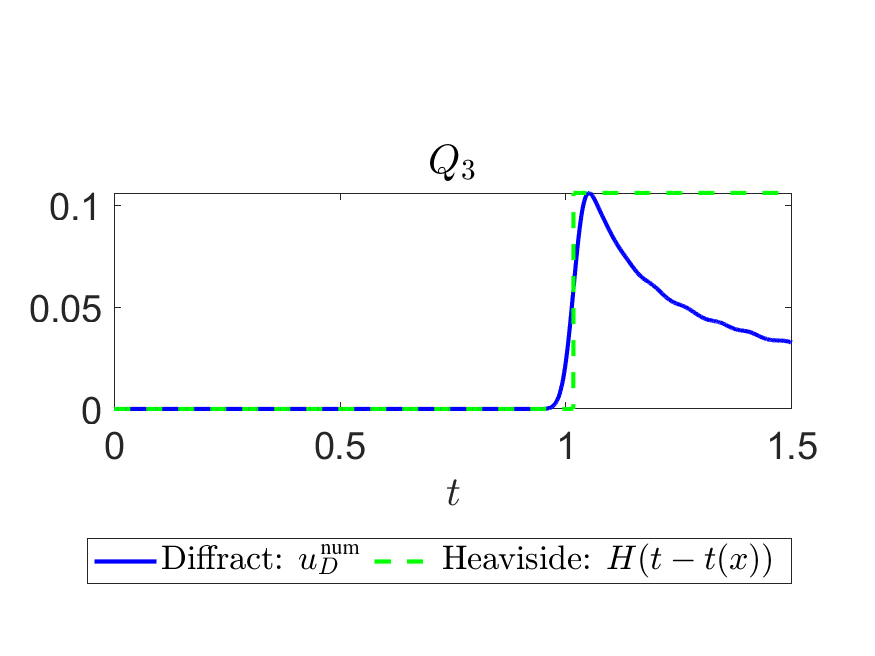}
\includegraphics[width=0.3\textwidth,trim=0 20 0 65,clip]{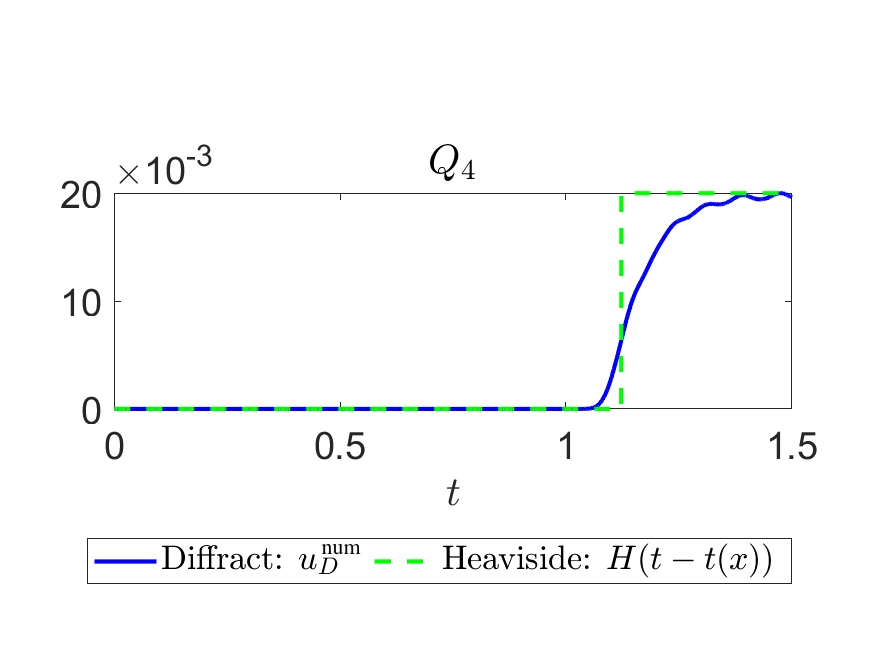}
\includegraphics[width=0.3\textwidth,trim=0 20 0 65,clip]{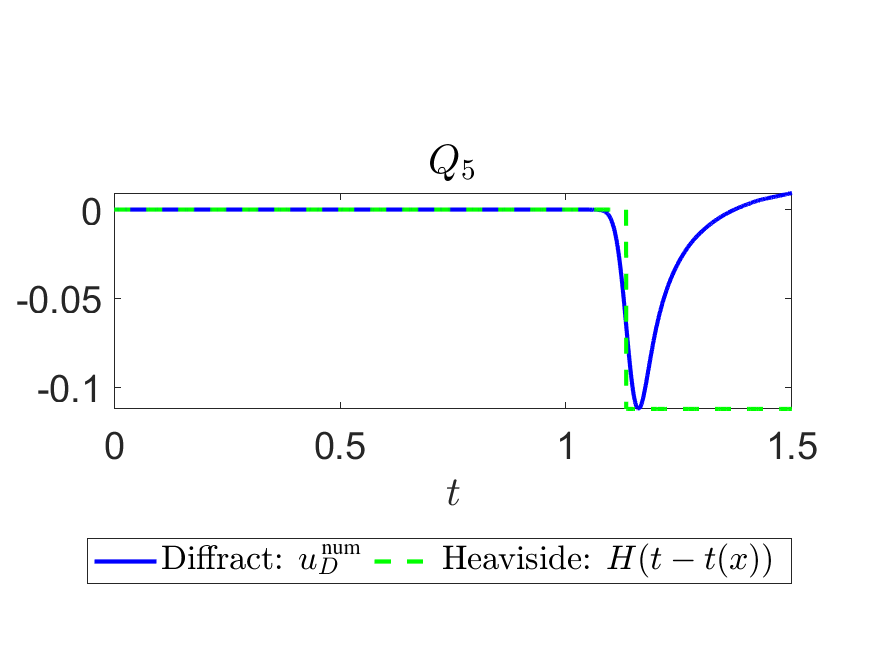}\\
\caption{Region B: positions $Q_{3}$, $Q_{4}$, and $Q_{5}$. %
After the incident signal, the diffracted signal dominates; the reflected signal is absent.}
\label{fig:signal-B}
\end{figure}

\begin{figure}[htb]
\centering
\includegraphics[width=0.32\textwidth,trim=0 20 0 65,clip]{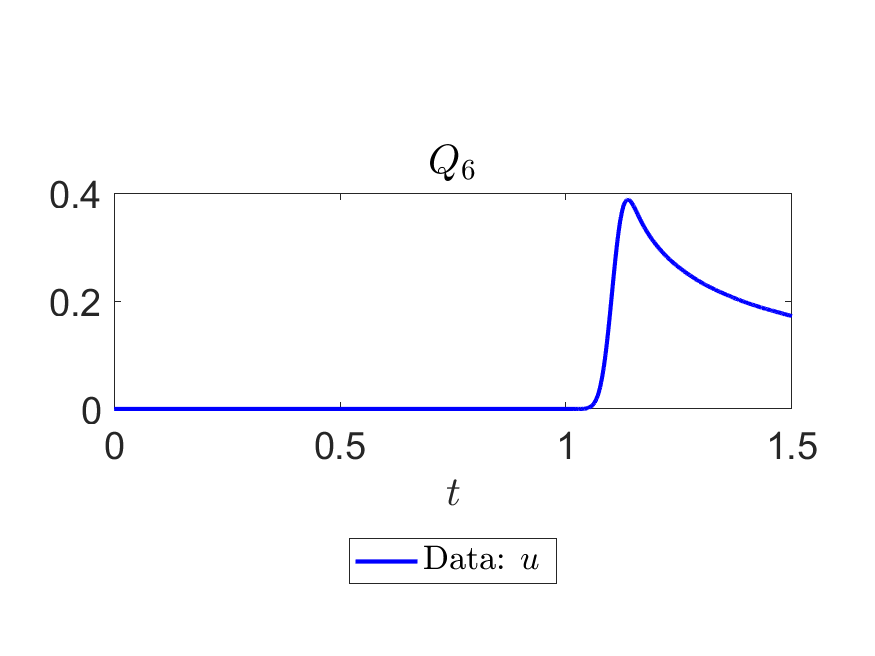}\\
\includegraphics[width=0.32\textwidth,trim=0 20 0 65,clip]{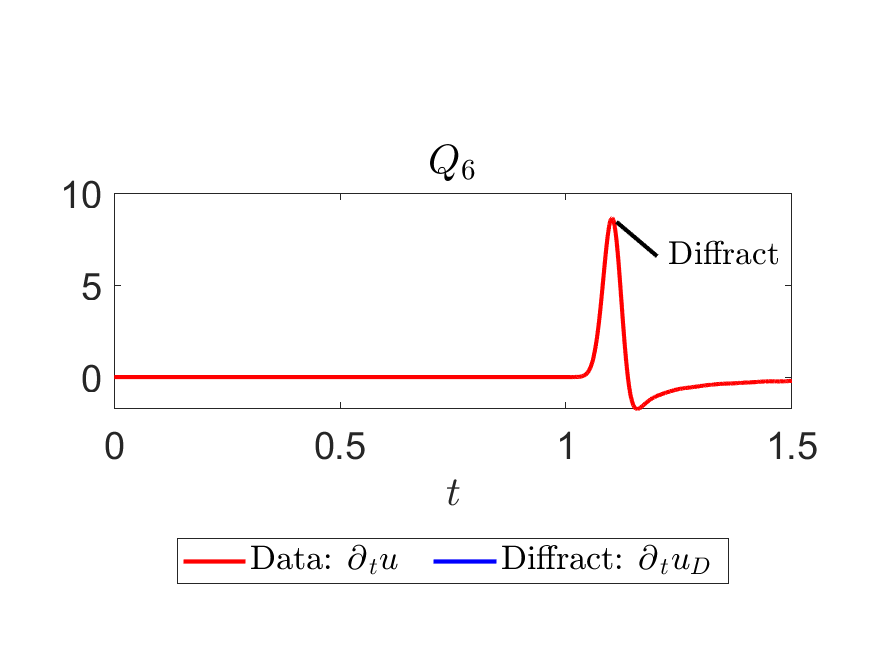}\\
\includegraphics[width=0.32\textwidth,trim=0 20 0 65,clip]{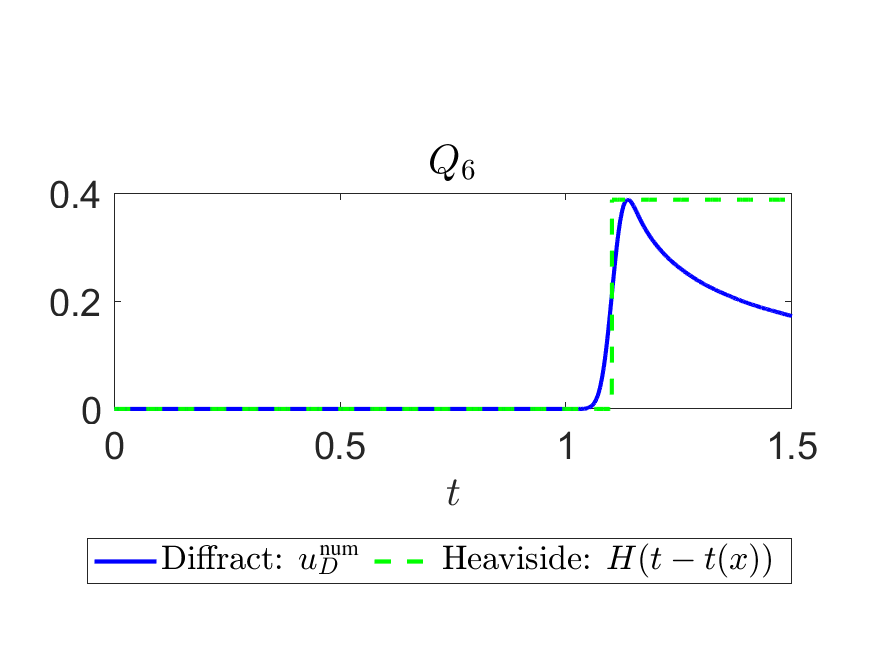}\\
\caption{Region C: position $Q_{6}$.
Only the diffracted signal is observed.}
\label{fig:signal-C}
\end{figure}

Signals at representative receivers in Regions A, B and C are shown in Figures \ref{fig:signal-A}--\ref{fig:signal-C}. %
In each panel we display the total wave field $u$, its time derivative, and the extracted diffracted signal $u_{\mathcal{D}}^{\mathrm{num}}$. %
In Region A (Figure \ref{fig:signal-A}, positions $Q_{1}$, $Q_{2}$) one typically observes three successive signals: the incident signal, a strong reflected signal, and a weaker diffracted signal. %
In Region B (Figure \ref{fig:signal-B}, positions $Q_{3}$, $Q_{4}$, $Q_{5}$) the reflected signal is missing; after the incident signal, the diffracted signal dominates and arrives at approximately $r + r_{1}(x)$. %
In Region C (Figure \ref{fig:signal-C}, position $Q_{6}$) the pattern is again diffraction-dominated, with no incident signal recorded. %
These observations confirm that diffracted signals remain available in angular sectors where reflected signals are absent. %

\subsection{Numerical inversion}
\label{subsec:inverse-num}

We now implement the reconstruction procedures of Theorems \ref{thm:uniqueness} and \ref{thm:polygon} on synthetic data. %
Each experiment is designed under the hypotheses of \cite{BCLX1}: a single pulsed source $O$, an exterior receiver curve $L$ that is not confined to hyperbolic branches, and reconstruction of only those tips that are visible. %
Unless otherwise stated, $L$ consists of one or several circular arcs (distinct from the full-circle receiver curve of Section \ref{subsec:diff-sim}, which serves only for visualisation of Regions A, B, C). %
The scattered field is obtained by subtracting the incident field as in \eqref{eqn:extract-diff}, and the primary diffracted signals are then isolated to yield $u_{\mathcal{D}}^{\mathrm{num}}$. %
For each diffracted signal, the arrival time $T_{i}$ is identified as the onset of the jump, detected via a threshold on $\lvert \partial_{t} u_{\mathcal{D}}^{\mathrm{num}} \rvert$ after the arrival of the geometric waves. %
The corresponding amplitude $U_{i}$ is then measured as the signed jump height. %
The numerical reconstruction follows Algorithm \ref{alg:two-step}. %
Both least-squares problems are solved by a trust-region reflective method with several initial guesses for the angular variables to reduce the effect of local minima. %
No wave equations are solved during the inversion. %

\begin{algorithm}[t]
\caption{Two-step reconstruction from diffracted signals}
\label{alg:two-step}
\begin{algorithmic}[1]
  \REQUIRE Source location $O$; receiver samples $\{ x_{i} \}_{i=1}^{M} \subset L$; numerical diffracted signals $u_{\mathcal{D}}^{\mathrm{num}}$ at $x_{i}$.
  \ENSURE Visible tip(s) $P$ and parameters $(\phi,\beta)$ (or $\{ P_{j},\phi_{j},\beta_{j} \}$ for several vertices).
  \STATE Extract diffracted arrival times $\{ T_{i} \}$ and amplitudes $\{ U_{i} \}$ from $u_{\mathcal{D}}^{\mathrm{num}}$.
  \STATE \textbf{Step 1 (location).} Solve
  \begin{equation*}
    \min_{P \in \mathbb{R}^{2}} \sum_{i=1}^{M} \bigl\lvert t_{P}(x_{i}) - T_{i} \bigr\rvert^{2}.
  \end{equation*}
  \STATE Convert $\{ x_{i} \}$ to polar coordinates $(r_{1}(x_{i}),\theta_{1}(x_{i}))$ centred at $P$; set $r = |O-P|$.
  \STATE \textbf{Step 2 (angles).} Solve
  \begin{equation*}
    \min_{(\phi,\beta)} \sum_{i=1}^{M} \Biggl\lvert -\frac{S(\theta_{1}(x_{i});\phi,\beta)}{16\pi \beta \sqrt{r\,r_{1}(x_{i})}} - U_{i} \Biggr\rvert^{2}
  \end{equation*}
  over the admissible set $\beta \in (0,1] \setminus \{N^{-1} : N \in \mathbb{Z}^{+} \}$.
  \STATE For polygons, isolate diffracted signals from each visible vertex by diffracted arrival time windows and repeat Steps 1--2.
  \RETURN Reconstructed visible tip(s) and angular parameters.
\end{algorithmic}
\end{algorithm}


\subsubsection{Convex single corner}
\label{subsec:convex-inv}

\begin{figure}[htb]
\centering
\includegraphics[width=0.28\textwidth,trim=40 0 70 0,clip]{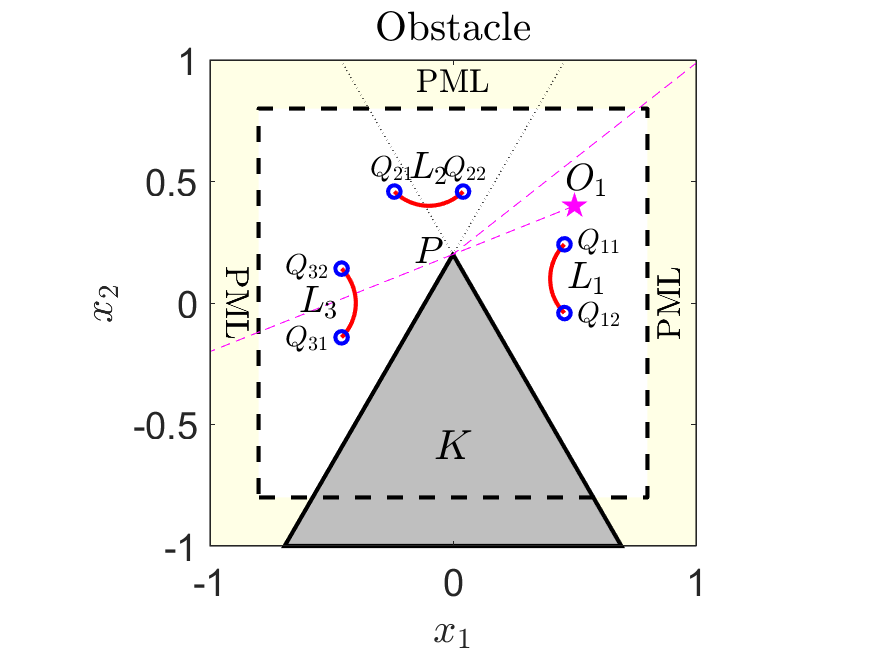}
\caption{Convex corner inversion geometry: source $O_{1}$ and arcs $L_{1}$, $L_{2}$, $L_{3}$.}
\label{fig:convex-inv-setup}
\end{figure}

We first reconstruct a single convex corner with true parameters $P^{\mathrm{true}} = (0,0.2)$ and $\beta^{\mathrm{true}} = \tfrac{5}{6}$, using the data measured on three exterior arcs $L_{1}$, $L_{2}$, $L_{3}$ (Figure \ref{fig:convex-inv-setup}). %

\begin{figure}[htb]
\centering
\includegraphics[width=0.2\textwidth,trim=100 0 100 0,clip]{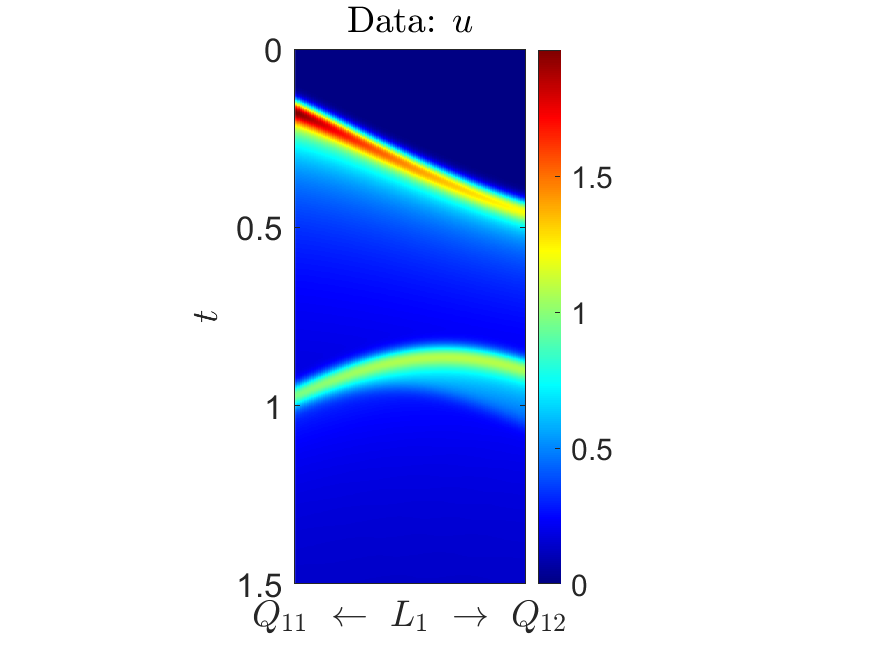}
\includegraphics[width=0.2\textwidth,trim=100 0 100 0,clip]{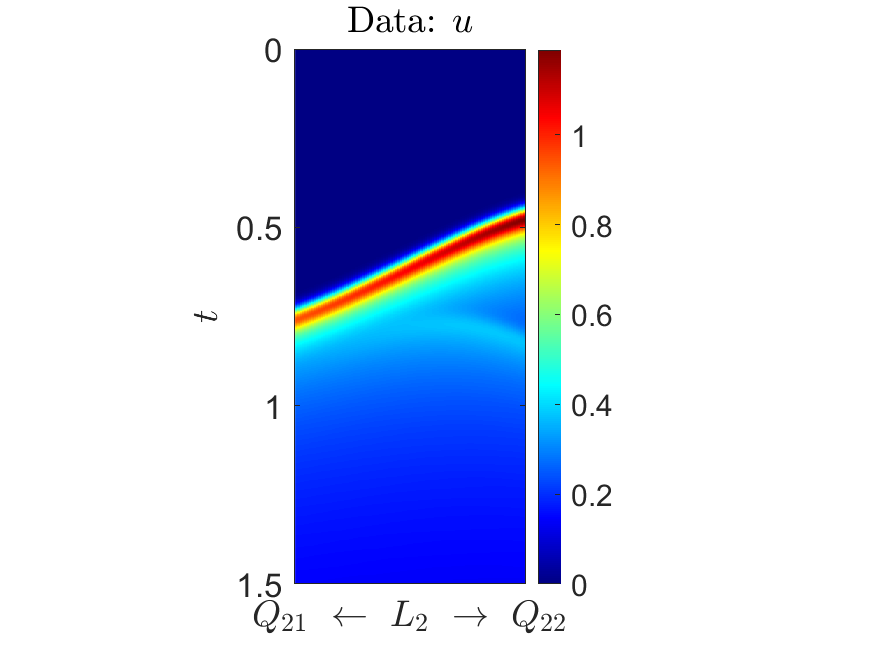}
\includegraphics[width=0.2\textwidth,trim=100 0 100 0,clip]{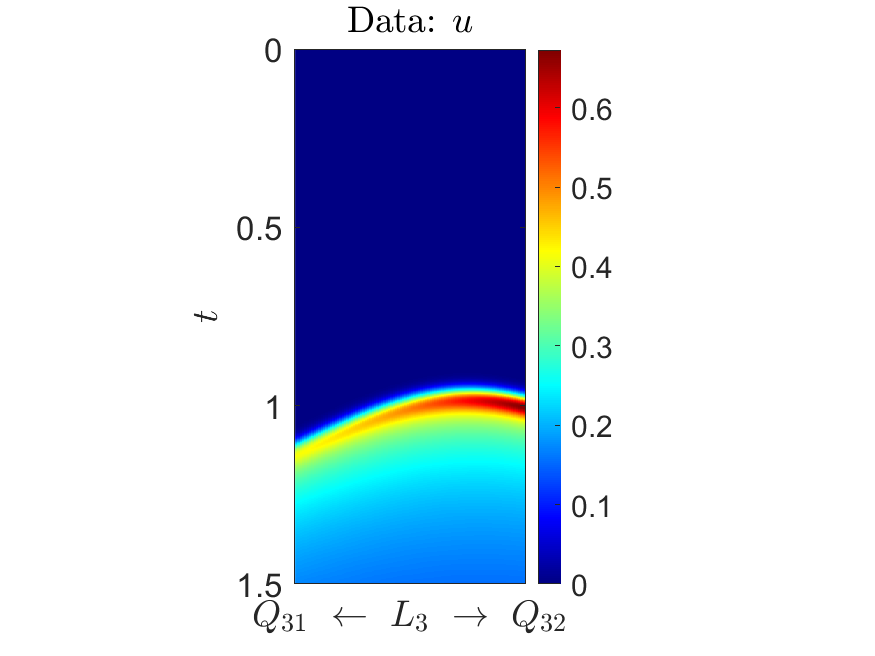}\\
\includegraphics[width=0.2\textwidth,trim=100 0 100 0,clip]{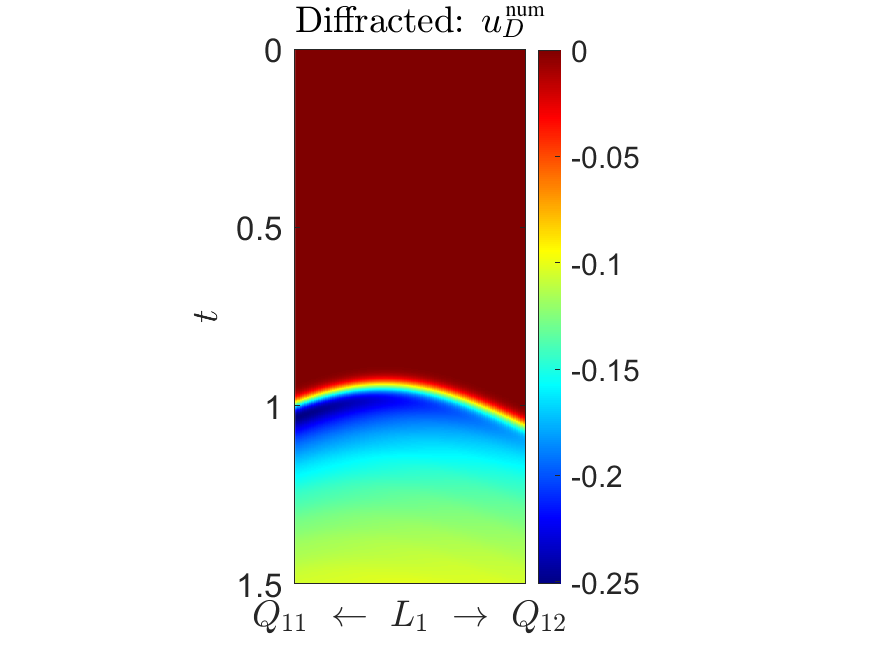}
\includegraphics[width=0.2\textwidth,trim=100 0 100 0,clip]{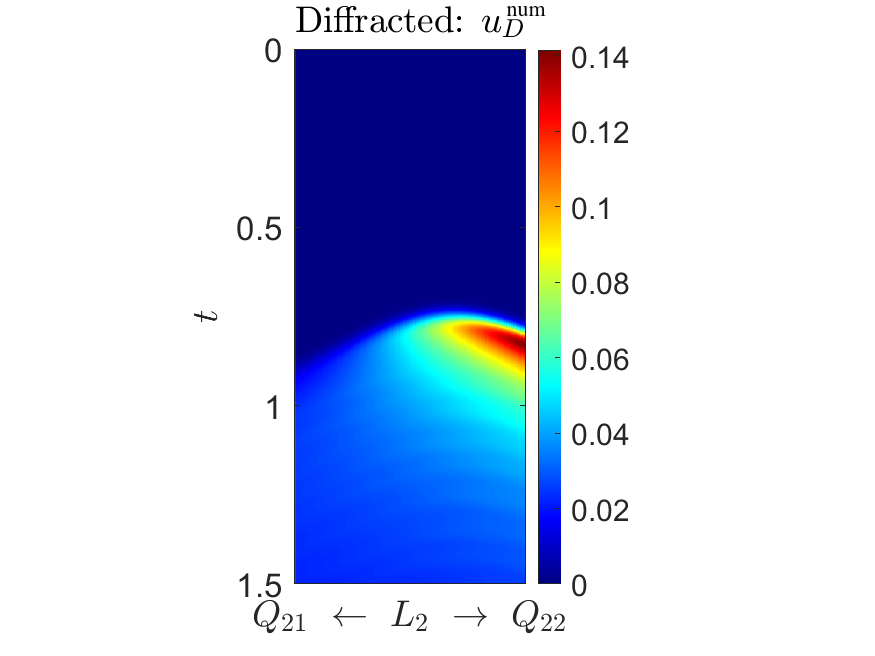}
\includegraphics[width=0.2\textwidth,trim=100 0 100 0,clip]{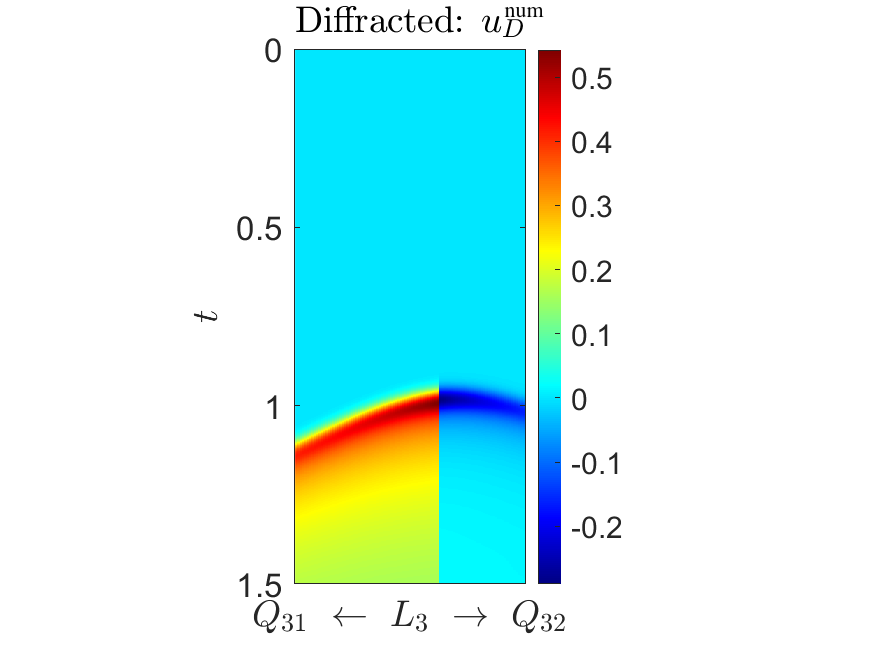}\\
\caption{Inversion data on $L_{1}$, $L_{2}$, $L_{3}$: total field (top) and numerical diffracted signals $u_{\mathcal{D}}^{\mathrm{num}}$ (bottom).}
\label{fig:convex-inv-traces}
\end{figure}

\begin{figure}[htb]
\centering
\includegraphics[width=0.32\textwidth,trim=10 0 30 0,clip]{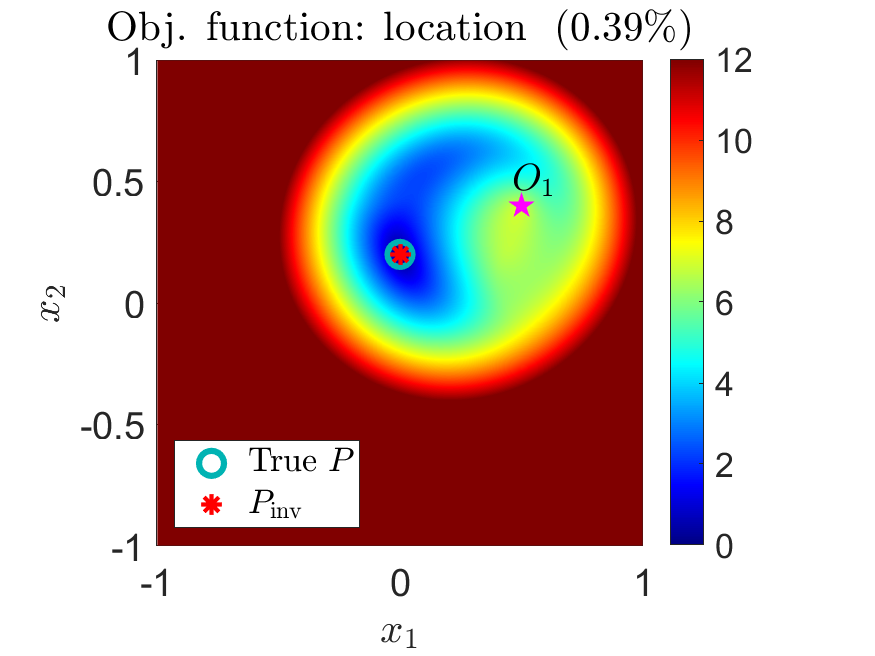}
\includegraphics[width=0.32\textwidth,trim=10 0 30 0,clip]{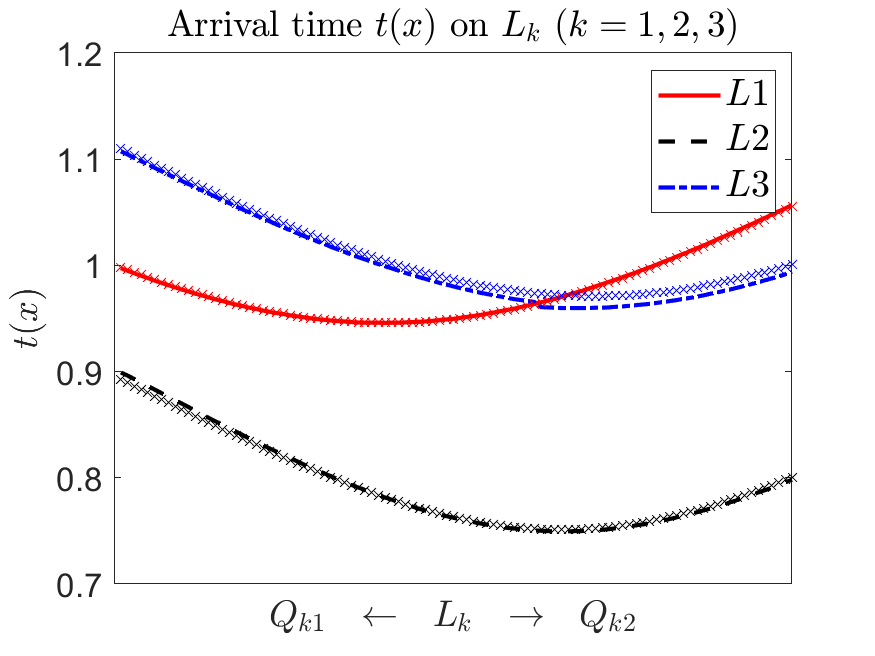}
\caption{Step 1 for the convex corner: arrival-time objective on $L_{1} \cup L_{2} \cup L_{3}$ (left) and arrival-time fit (right).}
\label{fig:convex-step1}
\end{figure}

Figure \ref{fig:convex-inv-traces} compares the total wave field on the three arcs with the extracted diffracted signals $u_{\mathcal{D}}^{\mathrm{num}}$ obtained from \eqref{eqn:extract-diff} after isolation. %
On each arc the diffracted signal is clearly identifiable after the arrival of geometric signals. %
Hence, we can extract both the arrival times and amplitudes of diffracted waves at the present resolution. %
The tip location is first estimated by solving the arrival-time least-squares problem \eqref{eqn:step1}. %
Figure \ref{fig:convex-step1} (left) shows a sharp minimum of the objective function near the true tip. %
The right panel shows that the predicted arrival times $\{ t_{P}(x_{i}) \}$ of diffracted waves match the picked times $\{ T_{i} \}$ closely. %

\begin{figure}[htb]
\centering
\includegraphics[width=0.28\textwidth,trim=40 0 70 0,clip]{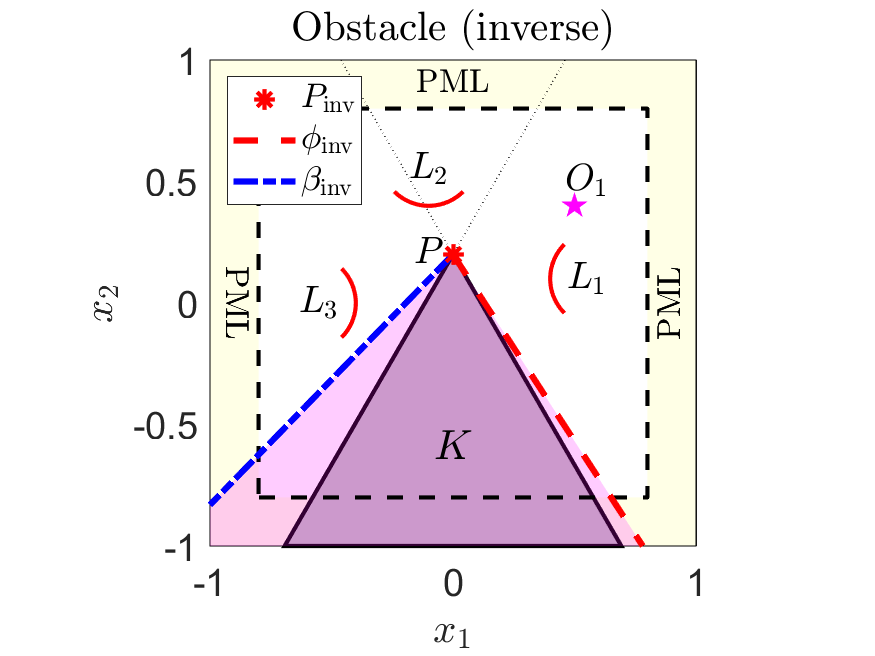}
\caption{Reconstructed convex corner after Steps 1--2 (arcs $L_{1} \cup L_{2} \cup L_{3}$).}
\label{fig:convex-recon}
\end{figure}

With $P$ fixed, the radial distances are known, and Step 2 fits $(\phi,\beta)$ through \eqref{eqn:step2}. %
The resulting wedge is plotted in Figure \ref{fig:convex-recon}. %
Quantitatively, the recovered tip has a relative error of about $0.39\%$ with respect to $P^{\mathrm{true}} = (0,0.2)$ in the computational window $[-1,1]^{2}$, while $\beta$ is recovered as approximately $0.825$ against $\beta^{\mathrm{true}} = 5/6 \approx 0.833$, i.e., a relative error of about one percent. %
These accuracies are consistent with the second-order spatial discretisation of the forward scheme and with the conditioning expected from the angular uniqueness result for $S(\theta_{1};\phi,\beta)$. %

\subsubsection{Convex corner with noisy data}
\label{subsec:convex-noise}

\begin{figure}[htb]
\centering
\includegraphics[width=0.2\textwidth,trim=100 0 100 0,clip]{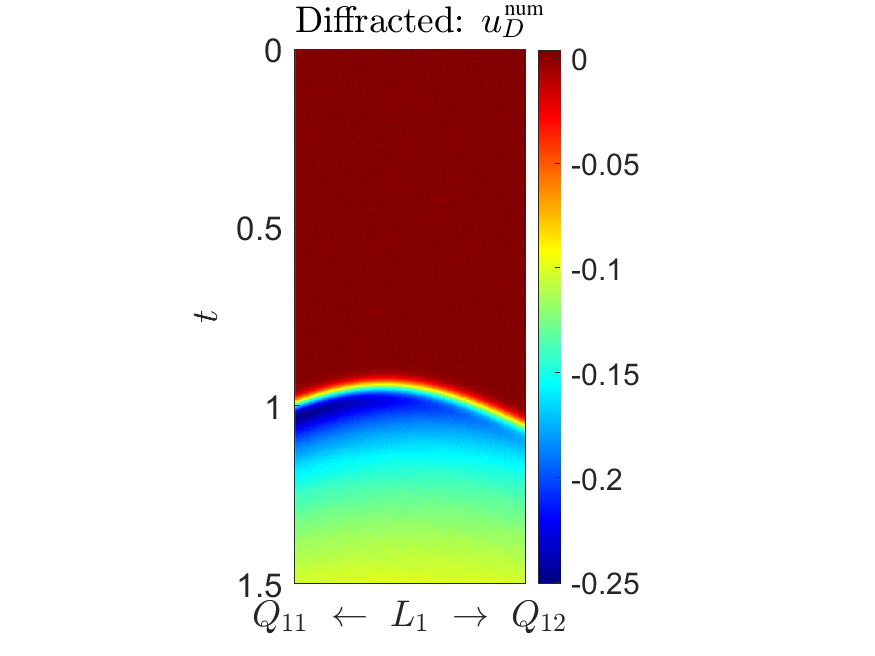}
\includegraphics[width=0.2\textwidth,trim=100 0 100 0,clip]{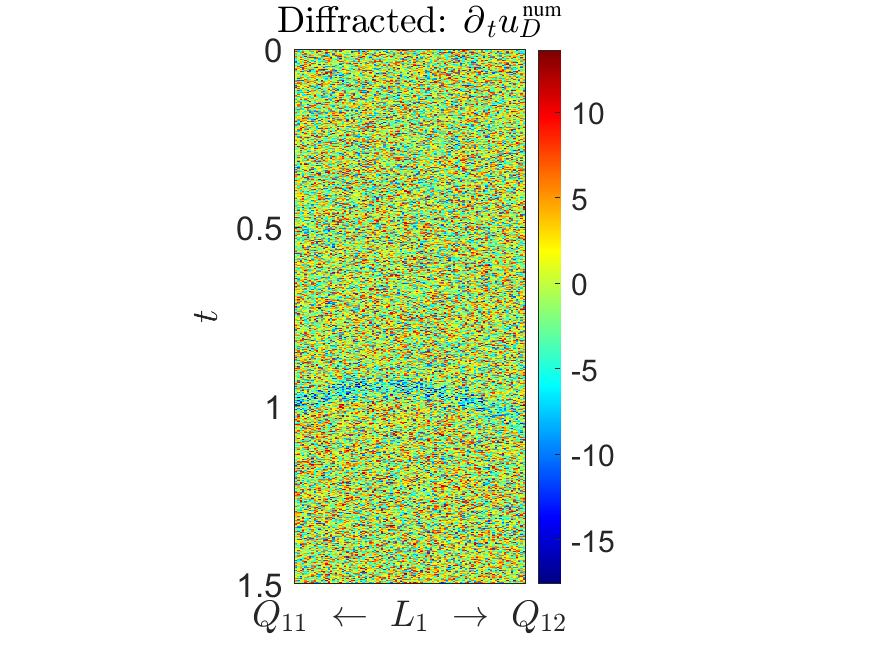}
\includegraphics[width=0.28\textwidth,trim=40 0 70 0,clip]{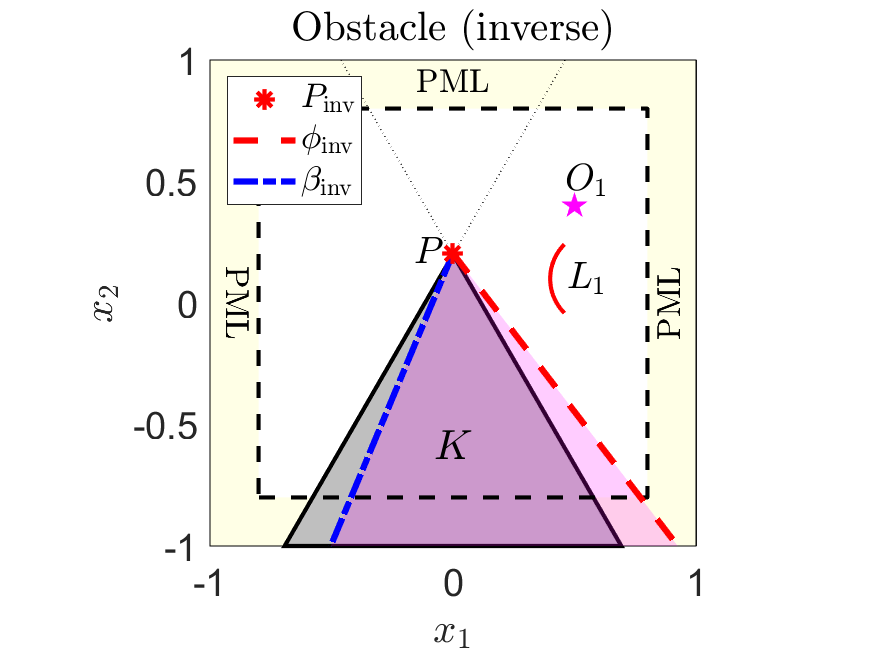}
\caption{Convex corner with noisy data on $L_{1}$ (noise level $0.2\%$). %
Left: numerical diffracted signal $u_{\mathcal{D}}^{\mathrm{num}}$. %
Centre: time derivative $\partial_{t} u_{\mathcal{D}}^{\mathrm{num}}$. %
Right: reconstructed corner after Steps 1--2.}
\label{fig:convex-noise}
\end{figure}

We next repeat the convex-corner experiment of Section \ref{subsec:convex-inv} after contaminating the extracted diffracted records by additive random noise at level $0.2\%$, restricting the inversion to the single exterior arc $L_{1}$. %
Figure \ref{fig:convex-noise} summarises the test: the noisy numerical diffracted signal $u_{\mathcal{D}}^{\mathrm{num}}$ and its time derivative on $L_{1}$, together with the reconstructed wedge after Algorithm \ref{alg:two-step}. %
The diffracted front remains identifiable, though the time derivative is more irregular than in the noise-free record of Figure \ref{fig:convex-inv-traces}. %
Step 1 on $L_{1}$ still produces a sharp minimum near the true tip. %
The relative tip error read from the objective plot is about $3.3\%$, compared with $0.39\%$ in the noise-free experiment of Section \ref{subsec:convex-inv}. %
The picked arrival times on $L_{1}$ remain close to the predicted travel times. %
With this tip fixed, Step 2 recovers $(\phi,\beta)$ with relative errors of about $9.7\%$ and $0.67\%$, respectively. %
Thus, at noise level $0.2\%$, the two-step procedure on the single arc $L_{1}$ continues to locate the tip and recover the opening. %

\subsubsection{Concave single corner}
\label{subsec:concave-num}

A concave corner corresponds to $\beta \in (0,\tfrac{1}{2})$. %
The formula for the principal part $u_{\mathcal{D}}$ in \eqref{eqn:diffract} remains valid, but the geometric content of Regions A, B, C changes. %
The illuminated and shadowed sectors are rearranged, and specular reflection may be absent on larger portions of the receiver curve. %
Provided the tip is diffractive and visible, the same inverse map as in Theorem \ref{thm:uniqueness} applies. %
Specifically, diffracted arrival times determine $P$, and amplitudes determine $(\phi,\beta)$. %

\begin{figure}[htb]
\centering
\includegraphics[width=0.28\textwidth,trim=40 0 70 0,clip]{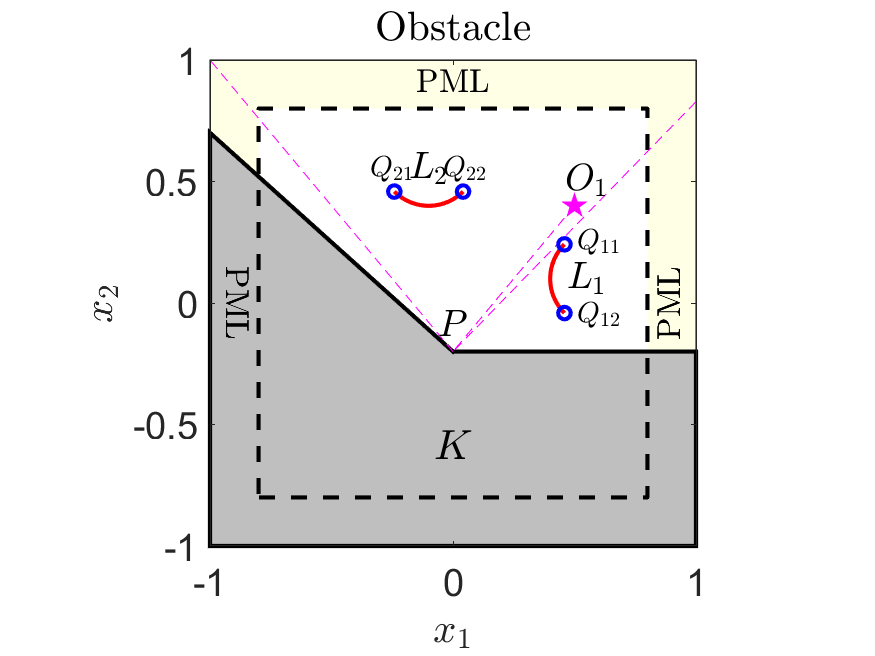}
\includegraphics[width=0.2\textwidth,trim=100 0 100 0,clip]{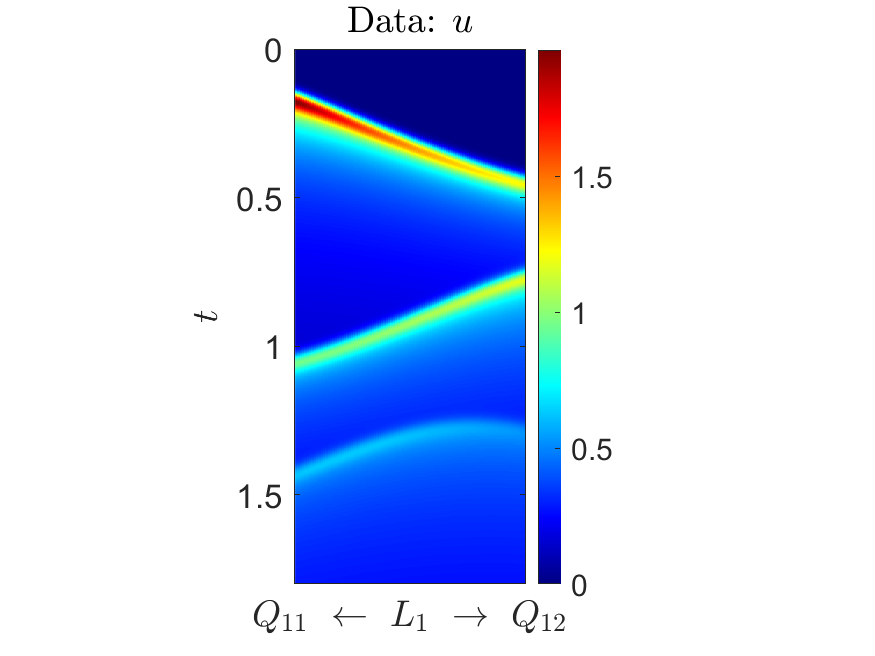}
\includegraphics[width=0.28\textwidth,trim=40 0 70 0,clip]{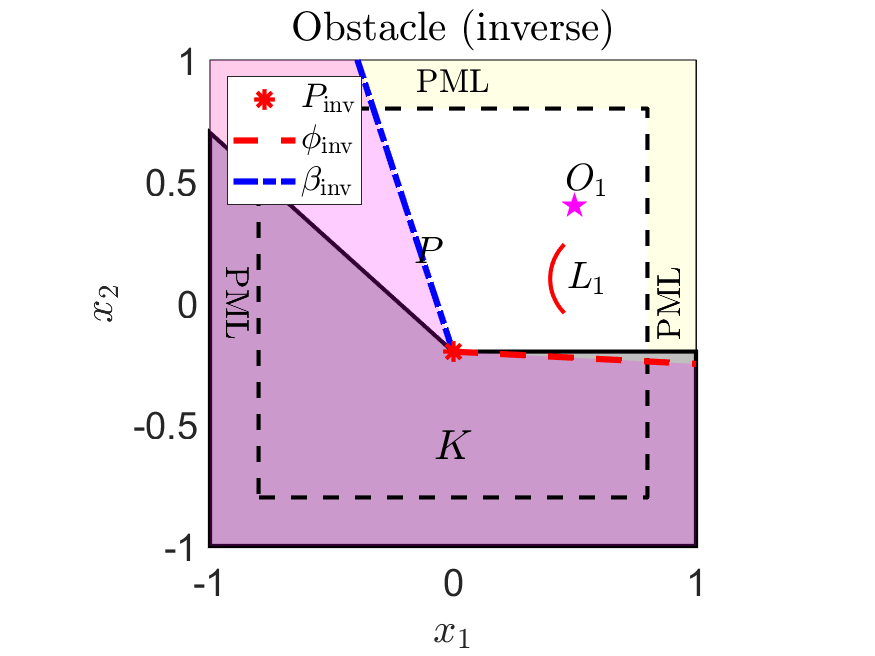}
\caption{Concave corner. %
Left: true geometry with source and receiver arc $L_{1}$. %
Centre: measured total field on $L_{1}$ (horizontal: receiver index; vertical: time). %
Right: reconstructed corner after Steps 1--2.}
\label{fig:concave}
\end{figure}

Figure \ref{fig:concave} summarises the experiment on a single concave tip. %
The left panel shows the obstacle, source and the arc $L_{1}$. %
The centre panel displays the total wave field recorded along $L_{1}$ as a function of receiver angle and time. %
The record is dominated by the incident and diffracted signals. %
The right panel shows the reconstructed corner after Algorithm \ref{alg:two-step}. %
The recovered tip and opening parameter $\beta$ match the true wedge within the resolution of the mesh, confirming that concavity alone does not obstruct the two-step procedure. %

\subsubsection{Triangular obstacle}
\label{subsec:triangle-num}

\begin{figure}[htb]
\centering
\includegraphics[width=0.28\textwidth,trim=40 0 70 0,clip]{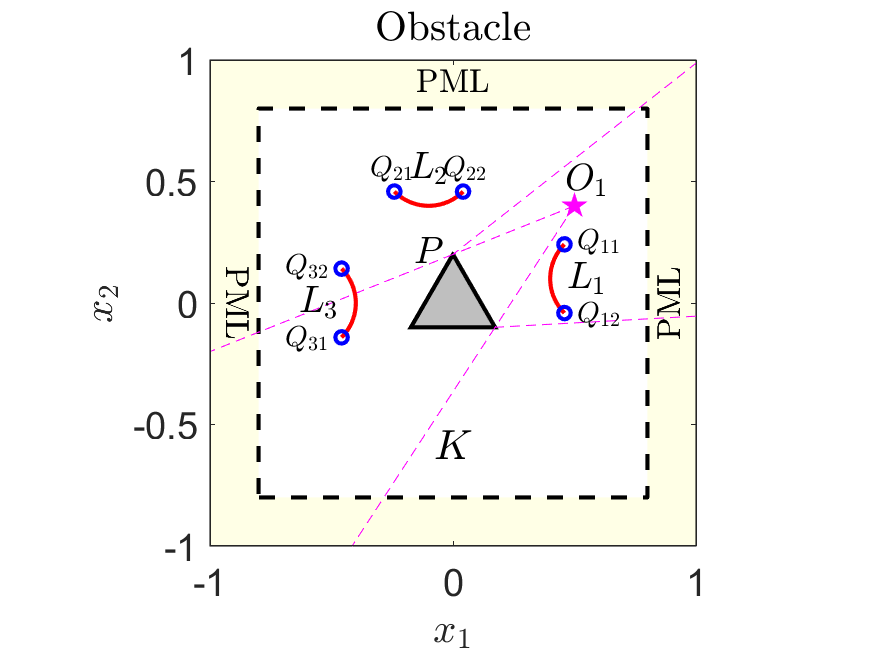}
\includegraphics[width=0.2\textwidth,trim=100 0 100 0,clip]{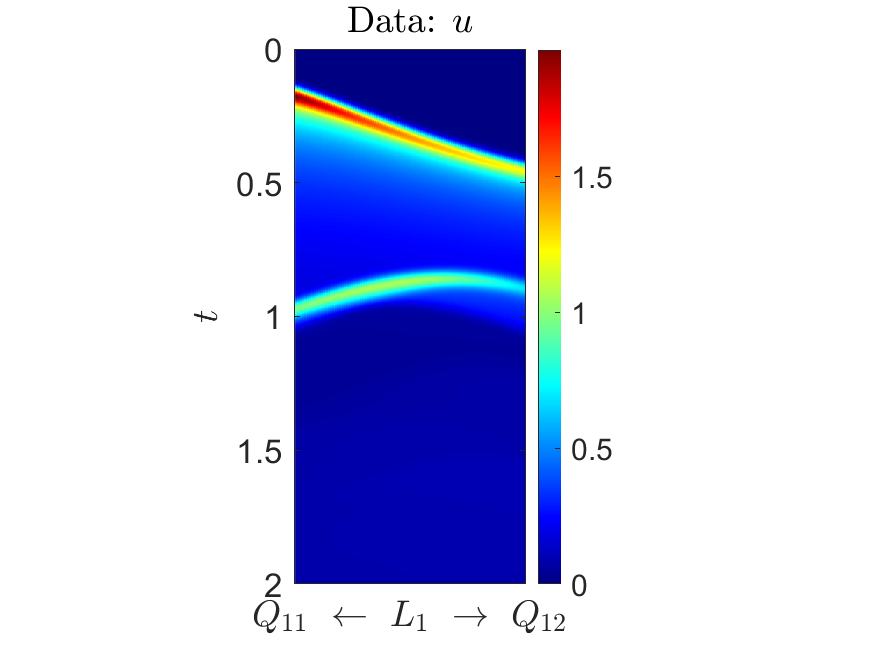}
\includegraphics[width=0.28\textwidth,trim=40 0 70 0,clip]{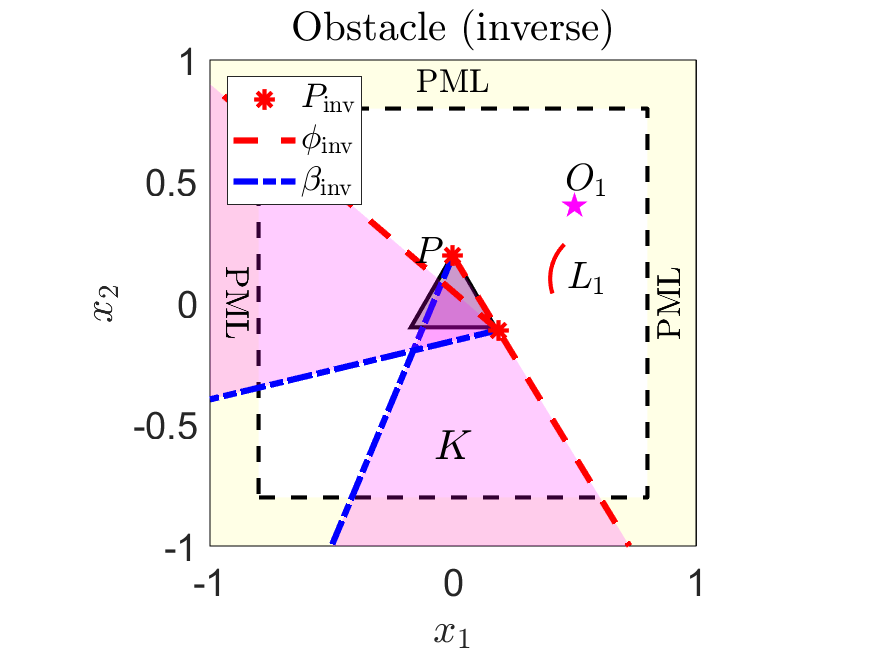}
\caption{Triangular obstacle. %
Left: true geometry with source $O_{1}$ and arcs $L_{1}$, $L_{2}$, $L_{3}$. %
Centre: total wave field on the inversion arc $L_{1}$. %
Right: reconstructed visible portion after sequential two-step inversion on $L_{1}$.}
\label{fig:triangle}
\end{figure}

We next consider a triangular obstacle. %
In the configuration of Figure \ref{fig:triangle} (left), the source $O_{1}$ and three exterior arcs $L_{1}$, $L_{2}$, $L_{3}$ are shown. %
The inversion reported below uses the primary diffracted signals recorded on $L_{1}$ to reconstruct the vertices visible with respect to $O_{1}$ and $L_{1}$. %
According to Section \ref{subsec:polygon-recovery}, each visible vertex generates a primary diffracted wave with its own diffracted arrival time $t_{j}(x)$. %
The space--time record on $L_{1}$ (Figure \ref{fig:triangle}, centre) therefore contains a sequence of diffracted signals arriving after the geometric signals. %
After isolating the primary diffracted signals and applying Algorithm \ref{alg:two-step} vertex by vertex, one obtains the reconstructed visible portion of the boundary displayed in Figure \ref{fig:triangle} (right). %

This experiment illustrates the reduction from Theorem \ref{thm:polygon} to computation. %
Specifically, the polygon is inverted not as a coupled multi-corner problem, but as a finite collection of independent single-corner problems whose solutions are then joined along the visible edges. %
Invisible vertices are not determined, in agreement with the theory. %

\subsubsection{Pentagonal obstacle}
\label{subsec:pentagon-num}

\begin{figure}[htb]
\centering
\includegraphics[width=0.28\textwidth,trim=120 0 120 0,clip]{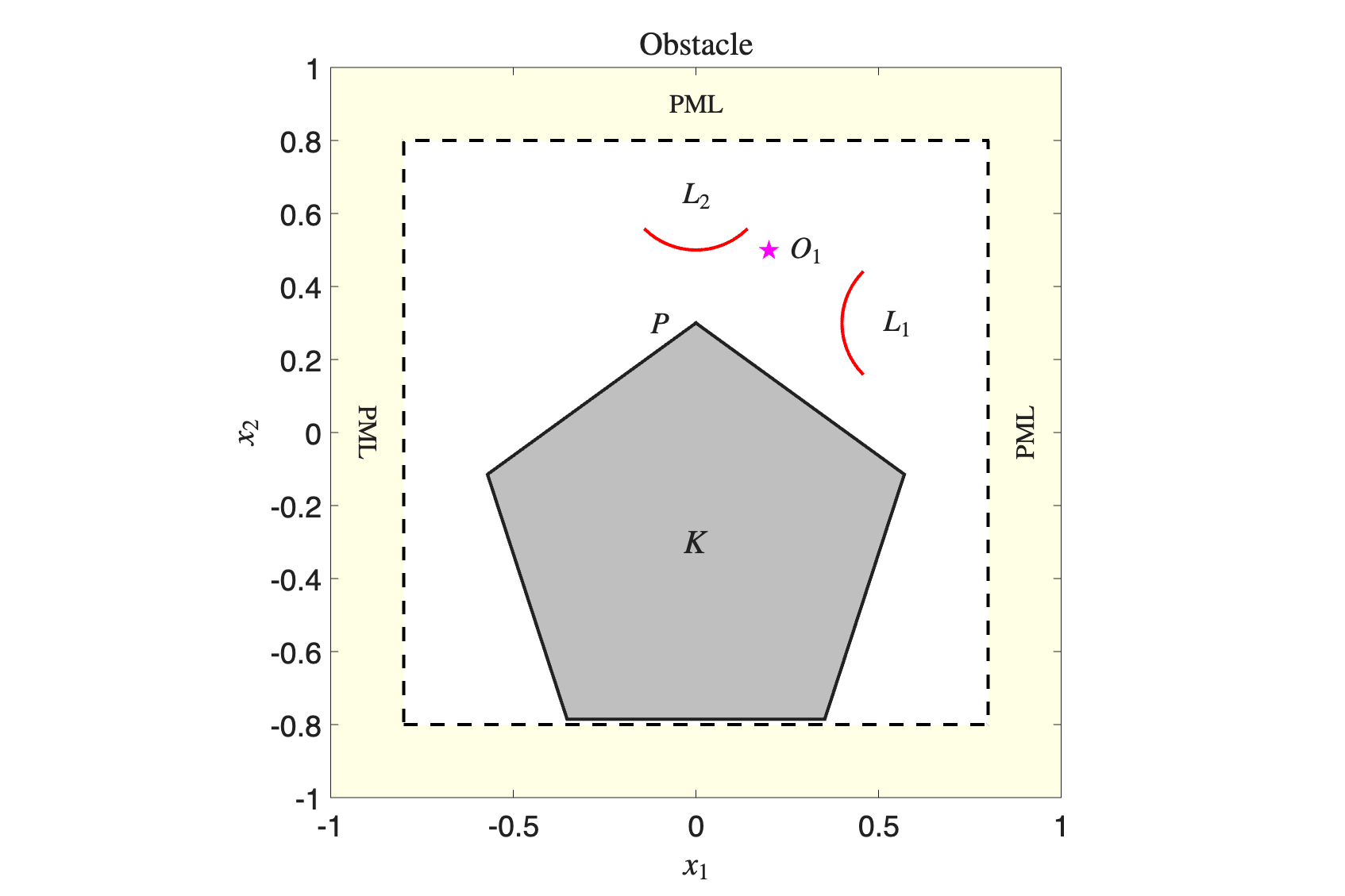}
\includegraphics[width=0.2\textwidth,trim=200 0 200 0,clip]{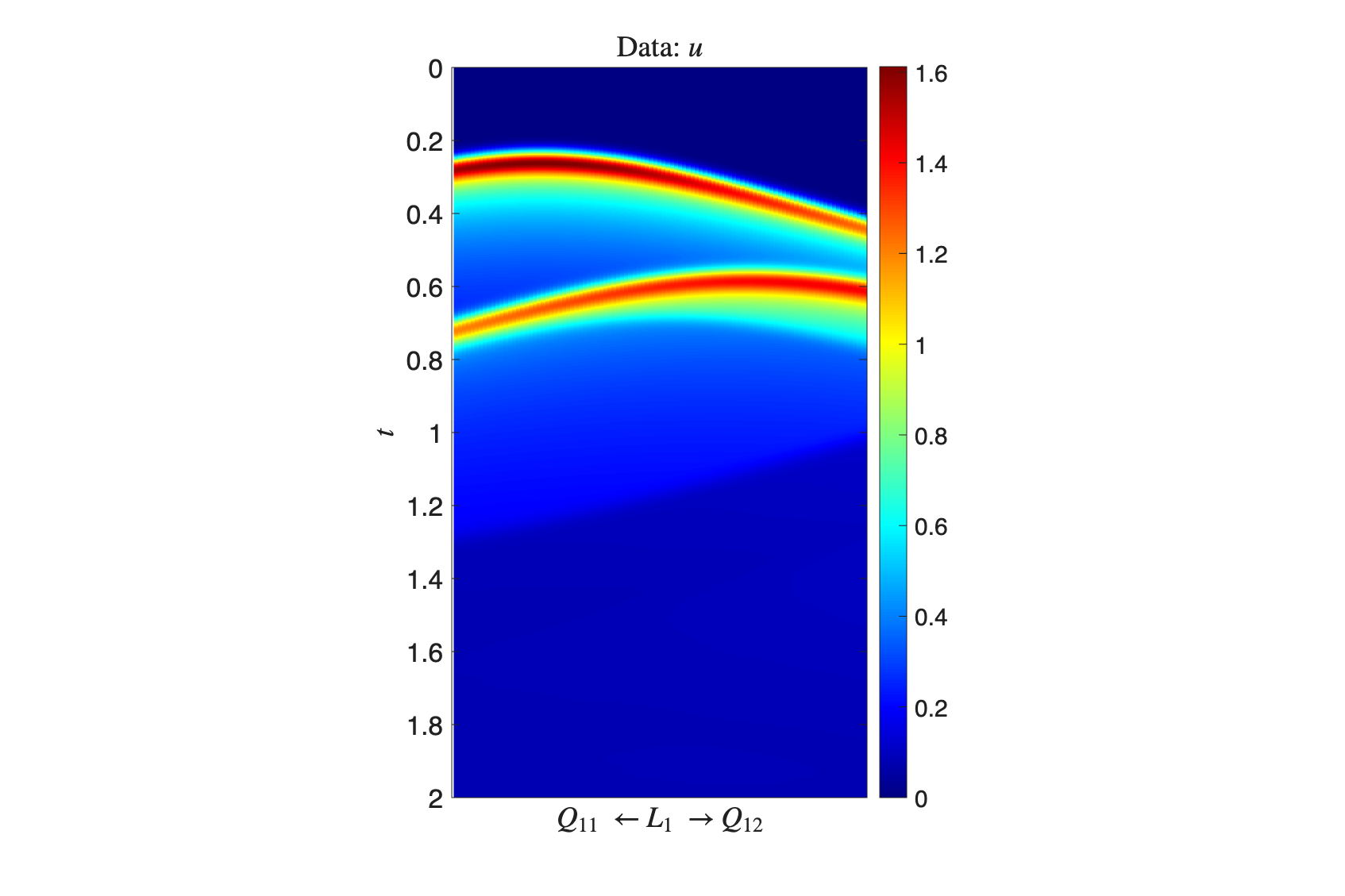}
\includegraphics[width=0.28\textwidth,trim=120 0 120 0,clip]{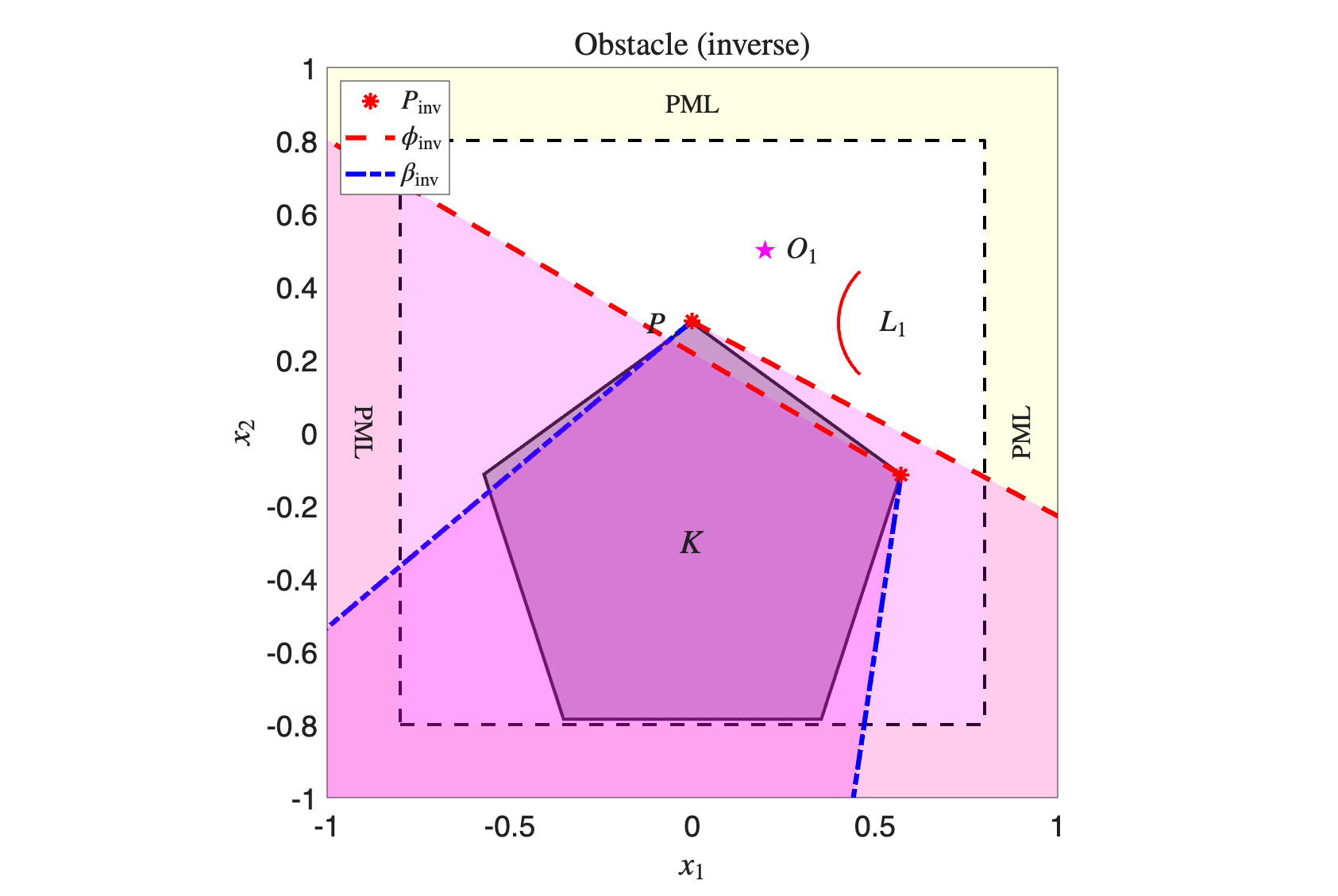}
\caption{Pentagonal obstacle. %
Left: true geometry with source $O_{1}$ and arcs $L_{1},L_{2}$. %
Centre: total field on the inversion arc $L_{1}$. %
Right: reconstructed visible vertices from $L_{1}$.}
\label{fig:pentagon}
\end{figure}

For a convex pentagon, the number of visible vertices typically increases, and the diffracted signals become more closely spaced in time. %
Figure \ref{fig:pentagon} (left) shows the true obstacle together with the source $O_{1}$ and two exterior arcs $L_{1}$, $L_{2}$. %
As in the triangular test, the inversion uses the $L_{1}$ record. %
The centre panel shows the total wave field recorded on $L_{1}$. %
The primary diffracted signals corresponding to visible vertices are extracted following the incident signal using arrival-time windows and the jump filter described above. %
We then reconstruct each visible corner separately and assemble the local wedges. %
The right panel of Figure \ref{fig:pentagon} displays the reconstructed visible vertices obtained from $L_{1}$ alone. %

This experiment leads to two observations. %
First, accurate picking of $\{ T_{i} \}$ remains possible because primary diffracted waves are more singular than multiple diffracted waves. %
Second, Step 2 is well conditioned on a single arc provided the range of receiver angles on $L_{1}$ is large enough for $S(\theta_{1};\phi,\beta)$ to vary non-trivially, as guaranteed by the injectivity result in Theorem \ref{thm:uniqueness}. %

\subsubsection{Star-shaped obstacle}
\label{subsec:star-num}

\begin{figure}[htb]
\centering
\includegraphics[width=0.28\textwidth,trim=120 0 120 0,clip]{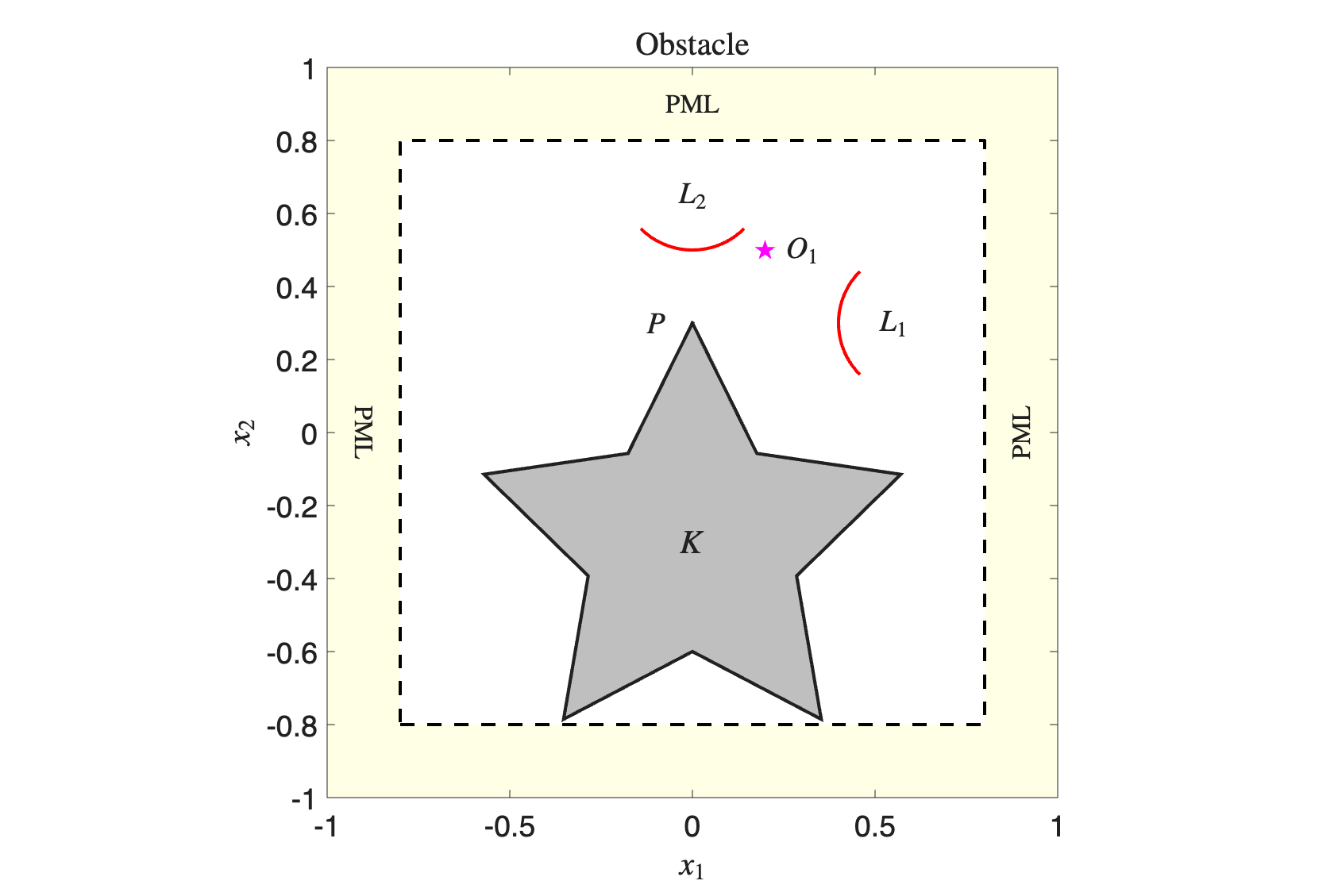}
\includegraphics[width=0.2\textwidth,trim=200 0 200 0,clip]{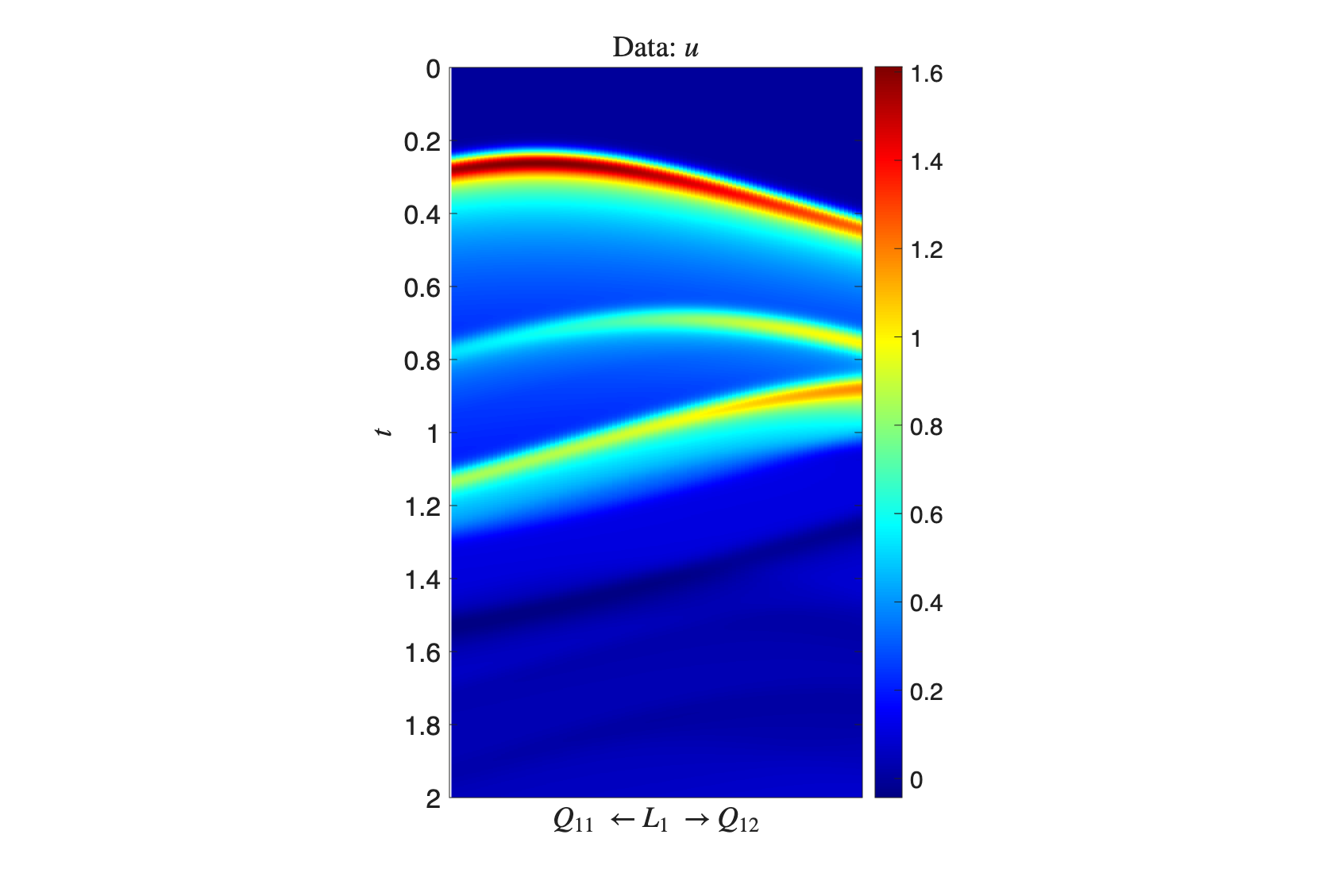}
\includegraphics[width=0.28\textwidth,trim=120 0 120 0,clip]{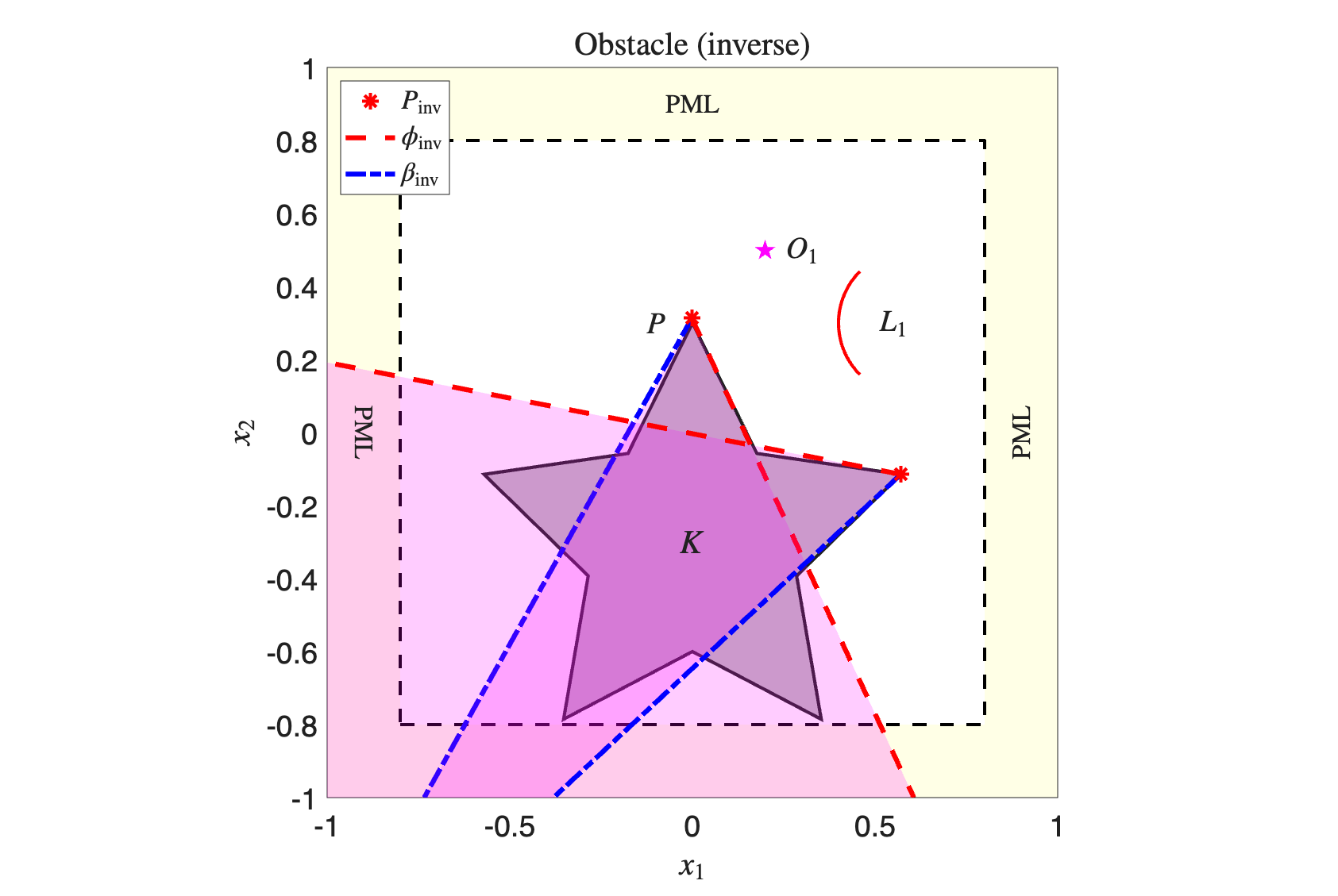}
\caption{Star-shaped obstacle. %
Left: true geometry with source $O_{1}$ and arcs $L_{1}$, $L_{2}$. %
Centre: total wave field on the inversion arc $L_{1}$. %
Right: reconstructed visible vertices from $L_{1}$.}
\label{fig:star}
\end{figure}

Finally we treat a non-convex star-shaped polygon, which mixes convex and concave corners. %
This geometry lies outside the convexity hypothesis of Theorem \ref{thm:polygon}, yet the local single-corner analysis of Theorem \ref{thm:uniqueness} still applies to every diffractive visible vertex. %
Figure \ref{fig:star} (left) shows the true star with source $O_{1}$ and arcs $L_{1}$ and $L_{2}$, and the inversion again uses $L_{1}$. %
The centre panel records the total wave field on $L_{1}$, in which the diffracted signals from several vertices appear successively. %
The right panel shows the reconstructed visible vertices after sequential application of Algorithm \ref{alg:two-step}. %

The experiment indicates that concavity at some vertices does not prevent reconstruction, as long as each reconstructed vertex remains visible and diffractive, and its diffracted signal can be separated in time. %


\section{Conclusions}
\label{sec:conclusions}

We have transformed the constructive diffraction-based theory into a computational reconstruction method for rigid polygonal obstacles. %
Diffracted signals are isolated from the total wave field and inverted in two nonlinear least-squares stages: arrival times of the diffracted waves locate visible vertices, and their amplitudes recover $(\phi,\beta)$ through the angular factor $S(\theta_{1};\phi,\beta)$. %
The angular partition into Regions A, B, C clarifies when specular reflection may be absent while diffraction persists, as described by the theory. %
Numerical experiments for a convex corner, a concave corner, a triangle, a pentagon, and a non-convex star-shaped obstacle support Theorems \ref{thm:uniqueness} and \ref{thm:polygon} employed here. %
Under additive noise at the level of $0.2\%$ on the same convex corner and the single arc $L_{1}$, the tip is recovered with a relative error of about $3.3\%$, and the opening parameter $\beta$ remains accurate to within one percent. %

Several limitations remain. %
The noise test uses a single convex corner, a single receiver arc and a known pulse shape. %
For vertices close to one another, arrival-time windows may overlap and a sharper microlocal filter would be required. %
Extensions to three-dimensional polyhedra, penetrable scatterers and multi-source experiments remain to be studied. %


\section*{Acknowledgments}

The AI tool was solely used to improve the English wording and readability of the manuscript. %
The authors independently carried out all mathematical reasoning and derivations, wrote and verified theorem statements and formulas, developed code, performed numerical experiments, and produced all figures. %
All AI-assisted language revisions have been reviewed by the authors. %
The authors take full responsibility for the entire content. %


\bibliographystyle{amsplain}
\bibliography{ref,numer}

\end{document}